\documentclass[10pt]{article}
\usepackage{amsmath,amsthm,amssymb,graphicx,xspace,epsfig,xcolor,colortbl,diagbox}

\usepackage{algorithm}\usepackage{algorithmic}
\usepackage{tikz, tkz-graph, tkz-berge}
\usetikzlibrary{decorations.pathreplacing}
\usetikzlibrary{patterns}
\usetikzlibrary{patterns.meta}
\usetikzlibrary{automata,arrows,positioning,calc}
\usepackage{pgf}
\usepackage{pgffor}
\usepackage{xspace}
\usepackage{color}
\usepackage{comment}
\usepackage{authblk}
\usepackage{multirow}
\usepackage{hyperref}
\usepackage{cleveref}
\usepackage{dsfont}
\usepackage{mathtools,amsfonts,amsthm,mathrsfs,amssymb,stmaryrd}
\usepackage{enumitem}
\usepackage{thm-restate}  
\usepackage{blkarray,booktabs,bigstrut}
\usepackage{subfigure}
\usepackage{hyperref}
\usepackage{listings}
\usepackage{xcolor}
\definecolor{codegreen}{rgb}{0,0.6,0}
\definecolor{codegray}{rgb}{0.5,0.5,0.5}
\definecolor{codepurple}{rgb}{0.58,0,0.82}
\definecolor{backcolour}{rgb}{0.95,0.95,0.92}
\lstdefinestyle{mystyle}{
    backgroundcolor=\color{backcolour},   
    commentstyle=\color{codegreen},
    keywordstyle=\color{blue},
    numberstyle=\tiny\color{codegray},
    stringstyle=\color{codepurple},
    basicstyle=\ttfamily\small,
    breakatwhitespace=false,         
    breaklines=true,                 
    captionpos=b,                    
    keepspaces=true,                 
    numbers=left,                    
    numbersep=5pt,                  
    showspaces=false,                
    showstringspaces=false,
    showtabs=false,                  
    tabsize=2
}
\usepackage{soul,cite}
\newtheorem{theorem}{Theorem}
\newtheorem{corollary}[theorem]{Corollary}
\newtheorem{proposition}[theorem]{Proposition}
\newtheorem{lemma}[theorem]{Lemma}

\theoremstyle{definition}

\newtheorem{conjecture}[theorem]{Conjecture}

\newtheorem{claim}[theorem]{Claim}

\newenvironment{proofclaim}[1][]
	{\par\noindent {\it Proof of the claim}. }{ \hfill$\lozenge$\par\vspace{11pt}}

\renewcommand{\epsilon}{\varepsilon}

\DeclareMathOperator{\UG}{UG}

\DeclareMathOperator{\OP}{OP}
\DeclareMathOperator{\op}{op}
\DeclareMathOperator{\dearth}{dearth}

\DeclareMathOperator{\dist}{dist}

\newcommand{\dic}{\vec{\chi}}
\newcommand{\bid}{\overleftrightarrow}
\newcommand{\Ra}{\Rightarrow}
\newcommand{\ind}[1]{[#1]}
\newcommand{\ori}[1]{\protect \overrightarrow{#1}}
\newcommand{\dis}{\boxplus}

\begin{document}

\title{The 3-dicritical semi-complete digraphs}
\author[,1]{Frédéric Havet\thanks{Research supported by the research grant DIGRAPHS ANR-19-CE48-0013}}
\author[,2]{Florian H{\"o}rsch\thanks{Research supported by the European Research Council (ERC) consolidator grant No.~725978 SYSTEMATICGRAPH.}}
\author[,1]{Lucas Picasarri-Arrieta\thanks{Research supported by the research grants DIGRAPHS ANR-19-CE48-0013 and ANR-17-EURE-0004.}}
\affil[1]{ Universit\'e C\^ote d'Azur, CNRS, Inria, I3S, Sophia-Antipolis, France}
\affil[2]{CISPA Saarbrücken, Germany}
\date{}
\maketitle
\vspace{-10mm}
\begin{center}
{\small 
\texttt{$\{$frederic.havet, lucas.picasarri-arrieta$\}$@inria.fr}\\ 
\texttt{florian.hoersch@cispa.de}\\
}
\end{center}

\maketitle

\begin{abstract}
A digraph is $3$-dicritical if it cannot be vertex-partitioned into two sets inducing acyclic digraphs, but each of its proper subdigraphs can. We give a human-readable proof that the number of 3-dicritical semi-complete digraphs is finite. Further, we give a computer-assisted proof of a full characterization of 3-dicritical semi-complete digraphs. There are eight such digraphs, two of which are tournaments. We finally give a general upper bound on the maximum number of arcs in a $3$-dicritical digraph.
\end{abstract}

\sloppy

\section{Introduction}

The purpose of this article is to completely characterize all semi-complete 3-dicritical digraphs, hence aiming to approach an analogue of some results on the maximum density of critical graphs. We first recall a few classical definitions and then give an overview of the previous work and motivation for our work.
\subsection{Definitions}

Our notation follows~\cite{bang2009}.
In this article, graphs and digraphs contain no parallel edges or arcs, respectively, and no loops. For some positive integer $k$, we use $[k]$ to denote the set $\{1,\ldots,k\}$.  The order of a graph $G$ (resp. digraph $D$) is denoted by $n(G)$ (resp. $n(D)$) and its number of arcs is denoted by $m(G)$ (resp. $m(D)$).

The {\it path} $v_1,\dots,v_n$ is the graph with vertex set $\{v_i \mid i\in [n]\}$ and edge set $\{v_iv_{i+1} \mid i\in [n-1]\}$. The {\it length} of a path is its number of edges. A {\it matching} is a set of pairwise disjoint edges.
A graph $H$ is a {\it subgraph} of a graph $G$ if $V(H)\subseteq V(G)$ and $E(H)\subseteq E(G)$ hold. 
If additionally $V(H)\neq V(G)$ or $E(H)\neq E(G)$, then $H$ is called a {\it proper} subgraph of $G$. 
For some $S \subseteq V(G)$, we define $G\ind{S}$ to be the subgraph of $G$ {\it induced by} $S$, that is, the graph whose vertex set is $S$ and whose edge set contains all the edges of $E(G)$ with both endvertices in $S$. 

Similarly, a digraph $D'$ is a {\it subdigraph} of a digraph $D$ if $V(D')\subseteq V(D)$ and $A(D')\subseteq A(D)$ hold. If additionally $V(D')\neq V(D)$ or $A(D')\neq A(D)$ holds, then $D'$ is called a {\it proper} subdigraph of $D$. If $V(D')=V(D)$, then $D$ is called a {\it spanning} subdigraph of $D$. 
For some $S \subseteq V(D)$, we define $D\ind{S}$ to be the subdigraph of $D$ {\it induced by} $S$, that is, the digraph whose vertex set is $S$ and whose arc set contains all the arcs of $A(D)$ with both endvertices in $S$. When we say that a digraph $D$ contains another digraph $D'$ as an (induced) subdigraph, we mean that $D$ has an (induced) subdigraph which is isomorphic to $D'$, hence not necessarily maintaining vertex labels. In case we want vertex labels to be maintained, we speak of an (induced) {\it labelled} subgraph. For two disjoint sets $X,Y \subseteq V(D)$, we say that $X$ {\it dominates} $Y$ (and $Y$ is dominated by $X$) if $xy \in A(D)$ holds for all $x \in X$ and $y \in Y$. 

The {\it underlying graph} of a digraph $D$, denoted by $\UG(D)$, is the graph on the same vertex set which contains an edge linking two vertices if the digraph contains at least one arc linking these two vertices.
Given an arc $uv$ in a digraph, we say that $v$ is an {\it out-neighbour} of $u$ and that $u$ is an {\it in-neighbour} of $v$. The set of out-neighbours of a vertex $u$ in a digraph $D$ is denoted by $N^+_D(u)$, and its set of in-neighbours in $D$ is denoted by $N^-_D(u)$. The out-degree of $u$ in $D$ is its number of out-neighbours, and its in-degree is its number of in-neighbours.

A \textit{digon} is a pair of arcs in opposite directions linking the same vertices. An arc is called {\it simple} if it is not contained in a digon.
A digraph in which every arc is contained in a digon is called {\it bidirected}.
An {\it oriented graph} is a digraph with no digon. A digraph $D$ is called {\it semi-complete}  if for all distinct $u,v \in V(D)$, at least one of the arcs $uv$ and $vu$ exists. A semi-complete digraph with no digon is called a {\it tournament}.
The {\it directed cycle} on $n\geq 2$ vertices is the digraph with vertex set $\{v_1,\dots,v_n\}$ and arc set $\{v_iv_{i+1} \mid i\in [n-1]\} \cup \{v_nv_1\}$. The directed cycle on three vertices is called the {\it directed triangle}.
A digraph is {\it acyclic} if it does not contain any directed cycle. For an acyclic digraph $D$, there is an ordering $v_1,\ldots,v_n$ of $V(D)$ such that for every arc $v_iv_j \in A(D)$, we have $i<j$. Such an ordering is called {\it acyclic}. An acyclic tournament is called {\it transitive}. For some positive integer $n$, we use $TT_n$ to denote the unique transitive tournament on $n$ vertices.

\subsection{Context}

A {\it $k$-colouring} of a graph $G$ is a function $\phi:V(G)\rightarrow [k]$. It is {\it proper} if $\phi(u)\neq \phi(v)$ for every edge $uv \in E(G)$. We say that $G$ is {\it $k$-colourable} if it admits a proper $k$-colouring. 
Colourability of graphs is one of the most deeply studied subjects in graph theory and countless aspects and variations of this parameter have been considered. One easy observation is that if a graph is $k$-colourable, then so is each of its subgraphs. 
This leads to the study of graphs that are in some way minimal obstructions to $(k-1)$-colourability, which are exactly $k$-critical graphs: a graph $G$ is {\it $k$-critical} if $G$ is not $(k-1)$-colourable, but all of its proper subgraphs are. Observe that every $k$-critical graph is $k$-colourable.
Dirac~\cite{Dirac52, Dirac55, Dirac57, Dirac67} established the basic properties of critical graphs.
For $k=1,2$, the only $k$-critical graph is $K_k$, the complete graph on $k$ vertices, and a graph is 3-critical if and only if it is an odd cycle. 
For every $k\geq 4$, there is no $k$-critical graph on $k+1$ vertices and for every $n \geq  k + 2$, there exists a $k$-critical graph on $n$ vertices.
Further, for every $k \geq 4$, the structure of $k$-critical graphs is immensely rich and many aspects of $k$-critical graphs have been studied, concerning both their importance in their own respect and their role when approaching colourability questions. 

One feature of particular interest is the density of $k$-critical graphs. For some $n \geq k$, let $g_k(n)$ be the minimum number of edges of a $k$-critical graph of order $n$ with the convention $g_k(k+1)=+\infty$.
The sequence $g_k$ is rather well understood due to the following result of Kostochka and Yancey~\cite{kostochkaJCTB109}.

\begin{theorem}[Kostochka and Yancey~\cite{kostochkaJCTB109}]\label{thm:g} Let $n$ and $k$ be two integers such that $n>k\geq 4$. 

If $G$ is a $k$-critical graph on $n$ vertices, then $m(G)\geq \left\lceil\frac{(k+1)(k-2)n-k(k-3)}{2(k-1)}\right\rceil$.
In other words, $g_k(n)\geq \left\lceil\frac{(k+1)(k-2)n-k(k-3)}{2(k-1)}\right\rceil$.
\end{theorem}
This lower bound is sharp for a significant amount of values of $n$ and $k$ and close to sharp for all remaining ones, which is certified by $k$-critical graphs that can be obtained from complete graphs using a construction due to  Haj\'os~\cite{hajos61}.
Theorem~\ref{thm:g} confirms a conjecture of Gallai~\cite{gallai63a} and improves on a collection of earlier results~\cite{Dirac57,gallai63b,gallai63a,krivelevich1997minimal,kostochka2003new}. A more detailed overview can be found in~\cite{kostochkaJCTB109}.  

\medskip

It is also interesting to determine the maximum number of edges in a $k$-critical graph. 
Erd\H{o}s~\cite{erdos69} asked, for every fixed $k\geq 4$, whether there exists a constant $c_k >0$ such that there exist arbitrarily large $k$-critical graphs $G$ with at least $c_k\cdot n(G)^2$ edges. This was proved by Dirac~\cite{Dirac52} when $k\geq 6$ and then by Toft~\cite{toftSSMH5} when $k\in \{4,5\}$. 
This initiated the quest after the supremum $c_k^*$, for fixed $k\geq 4$, of all values $c_k$ for which the statement holds. The following lower bound on $c_k^*$ follows from the explicit construction given in~\cite{toftSSMH5} and is still the best current bound.

\begin{theorem}[Toft~\cite{toftSSMH5}]
	For every integer $k\geq 4$ and infinitely many values of $n$, there exists a $k$-critical graph with $n$ vertices and at least $\frac{1}{2}\left(1 - \frac{3}{k-\delta_k} \right)n^2$ edges, where $\delta_k = 0$ if $k=0 \mod 3$, $\delta_k = \frac{4}{7}$ if $k=1 \mod 3$, and $\delta_k = \frac{22}{23}$ if $k=2 \mod 3$.
\end{theorem}

Concerning the upper bounds on $c_k^*$, observe that a $k$-critical graph does not contain any copy of $K_k$ as a subgraph. A seminal result of Tur\'an~\cite{turan1941} implies that such a graph $G$ of order $n$ has at most $\frac{1}{2}\left(1 - \frac{1}{k-1}\right)n^2$ edges (when $n = 0 \mod k$). Hence we have $c_k^* \leq \frac{1}{2}\left(1 - \frac{1}{k-1}\right)$. In 1987, Stiebitz~\cite{stiebitzComb7} improved on this lower bound.

\begin{theorem}[Stiebitz~\cite{stiebitzComb7}]
	For every integer $k\geq 4$ and sufficiently large integers $n$, every $k$-critical graph $G$ of order $n$ has at most $\frac{1}{2}\left(1 - \frac{1}{k-2}\right)n^2$ edges.
\end{theorem}

This remained the best upper bound on $c_k^*$ for many years, until Luo, Ma, and Yang~\cite{luoCPC32} qualitatively improved it in 2023.

\begin{theorem}[Luo, Ma, and Yang~\cite{luoCPC32}]
	For every integer $k\geq 4$, there exists $\epsilon_k \geq \frac{1}{18(k-1)^2}$ such that for sufficiently large integers $n$, every $k$-critical graph $G$ of order $n$ has at most $\frac{1}{2}\left(1 - \frac{1}{k-2} - \epsilon_k \right)n^2$ edges.
\end{theorem}

It remains an open problem to find the exact value of $c_k^*$. However, when $k\geq 6$, the analogue of $c_k^*$ is well-understood for triangle-free graphs. Indeed, for $k\geq 6$, Pegden~\cite{pegdenComb33} proved that there exist infinitely many $k$-critical triangle-free graphs $G$ with $\left(\frac{1}{4} - o(1)\right)n(G)^2$ edges. This is asymptotically best possible because of Turan's result.

\medskip

Several analogues of colouring have been introduced for digraphs. In 1982, Neumann-Lara~\cite{neumannlaraJCTB33} introduced the one of dicolouring.
A {\it $k$-dicolouring} of a digraph $D$ is a function $\phi:V(D)\rightarrow [k]$ such that $D\ind{\phi^{-1}(i)}$ is acyclic for every $i \in [k]$. We say that $D$ is {\it $k$-dicolourable} if it admits a $k$-dicolouring. 
Dicolourablity is a generalization of colouring to digraphs: indeed there is a trivial one-to-one correspondence between the proper $k$-colourings of a graph $G$ and the $k$-dicolourings of the associated bidirected graph $\bid{G}$ obtained from $G$ by replacing every edge by a digon.

We say that $D$ is {\it $k$-dicritical} if $D$ is not $(k-1)$-dicolourable, but all of its proper subdigraphs are. Observe that every $k$-dicritical digraph is $k$-dicolourable. 
The interest in $k$-dicritical graphs arises in a similar way as the interest in $k$-critical graphs. 
While the only 1-dicritical digraph is the digraph on one vertex and a graph is 2-dicritical if and only if it is a directed cycle, already 3-dicritical digraphs have a very diverse structure. 
Analogues of Haj\'os' construction have been found by Bang-Jensen et al.~\cite{bangjensenBSS2019}. Again, for some $n \geq k$, it is natural to consider $d_k(n)$, the minimum number of arcs of a $k$-dicritical digraph of order $n$, with the convention $d_k(n)=+\infty$ if no such digraph exists.

Observe that, for every $k\geq 2$, the digraph obtained from a bidirected complete graph on $k+1$ vertices by removing the arcs of a directed triangle is a $k$-dicritical digraph on $k+1$ vertices. 
Note also that a graph $G$ is $k$-critical if and only if its associated bidirected graph $\bid{G}$ is $k$-dicritical.
Together with the fact that there exist $k$-critical graphs of every order at least $k+2$,  this implies that there exists a $k$-dicritical digraph of any order at least $k$. In other words, we have $d_k(n) < +\infty$ for all $n\geq k$. Moreover, it holds $d_k(n)\leq 2g_k(n)$ for all $n\geq k$.

Kostochka and Stiebitz~\cite{kostochkaGC36} conjectured that the bidirected $k$-dicritical digraphs obtained from $k$-critical graphs are the sparsest $k$-dicritical digraphs.

\begin{conjecture}[Kostochka and Stiebitz~\cite{kostochkaGC36}]\label{ksconj1}
Let $n$ and $k$ be two integers such that $n-2\geq k\geq 4$.
Then $d_k(n) = 2g_k(n)$ and the $k$-dicritical digraphs of order $n$ with $d_k(n)$ arcs are the bidirected graphs associated to $k$-critical graphs with $g_k(n)$ edges. 
\end{conjecture}

With Theorem~\ref{thm:g} this conjecture implies the following sligthtly weaker one.

\begin{conjecture}[Kostochka and Stiebitz~\cite{kostochkaGC36}]\label{ksconj2}
Let $n$ and $k$ be two integers such that $n>k\geq 4$. If $D$ is a $k$-dicritical digraph on $n$ vertices, then $m(D)\geq \left\lceil\frac{(k+1)(k-2)n-k(k-3)}{k-1}\right\rceil$ and all digraphs attaining this bound are bidirected.
\end{conjecture}

In~\cite{kostochkaGC36}, Kostochka and Stiebitz confirmed Conjecture~\ref{ksconj2} for $k \leq 4$.

\medskip

It is expected that the  minimum number of arcs in a $k$-dicritical digraph
of order $n$ is larger than $d_k(n)$ if we impose this digraph to have no short directed cycles, and in particular if the digraph is an oriented graph.
Let $o_k(n)$ denote the minimum number of arcs in a $k$-dicritical oriented graph of order $n$ with the convention $o_k(n)=+\infty$ if there is no $k$-dicritical
oriented graph of order $n$. 
Clearly, we have $o_k(n) \geq d_k(n)$.
Kostochka and Stiebitz~\cite{kostochkaGC36} posed a conjecture suggesting that there is a significant gap between $d_k(n)$ and $o_k(n)$. 

\begin{conjecture}[Kostochka and Stiebitz~\cite{kostochkaGC36}]\label{ksconj3}
There exists $\epsilon>0$ such that
$o_k(n) \geq (1+\epsilon)\cdot d_k(n)$ for every $k\geq 4$ and sufficiently large $n$.
\end{conjecture}

Observe that this conjecture trivially holds when
$o_k(n)=+\infty$. However we do not know the set $N_k$ of integers $n$ for which $o_k(n) < +\infty$, that is, for which there exists a $k$-dicritical oriented graph on $n$ vertices. Already the minimum number $n_k$ of vertices of a $k$-dicritical oriented graph is unknown except for small values of $k$ :
clearly $n_2=3$ ; Neumann Lara~\cite{neumannlaraDM135} proved $n_3=7$ and $n_4=11$ ; 
Bellitto et al.~\cite{bellittoMC93} recently established $n_5=19$. 
As observed by Aboulker et al.~\cite{bellittoArXiv22} using a lemma of Hoshino and Kawarabayashi~\cite{hoshinoComb35}, there exists a smallest integer $p_k$ such that there exists a $k$-dicritical oriented graph on $n$ vertices for any $n\geq p_k$.
Moreover, while $p_3 = n_3=7$, they showed that $p_4 \neq n_4$ because there is no $4$-dicritical oriented graph on $12$ vertices.

Conjecture~\ref{ksconj3} has been confirmed for $k=3$ by Aboulker et al.~\cite{bellittoArXiv22} and for $k=4$ by the first and third authors, and Rambaud~\cite{havetArXiv23}.

\medskip
We now turn our attention to the maximum density of $k$-dicritical digraphs which is the main subject of the present article. 
For every $k \geq 3$, Hoshino and Kawarabayashi~\cite{hoshinoComb35} constructed an infinite family of $k$-dicritical oriented  graphs  $D$ on $n$ vertices  which satisfy $m(D)\geq (\frac{1}{2}-\frac{1}{2^k-1})n^2$, and they conjectured that this bound is tight. 

\begin{conjecture}[Hoshino and Kawarabayashi~\cite{hoshinoComb35}] 
\label{conj:hoshino_kawarabayashi}
Let $k\geq 3$ be an integer.
  If $D$ is a $k$-dicritical oriented graph, then $m(D)\leq (\frac{1}{2}-\frac{1}{2^k-1})n^2$.
\end{conjecture}

  Aboulker~\cite{perso_aboulker} observed that, since a tournament has $\frac{1}{2}n(n-1)$ arcs, this conjecture implies that the number of $k$-dicritical tournaments is finite, and he asked whether this latter statement holds.
 It trivially does for the case $k=2$.
 
 \subsection{Our results}
 
 In this paper, we positively answer  Aboulker's question in the case $k=3$ by showing that the collection of $3$-dicritical semi-complete digraphs is finite, and hence so is the subcollection of $3$-dicritical tournaments.

\begin{restatable}{theorem}{thmfinitethreedicriticalsemicomplete}
\label{thm:finite_3_dicritical_semi_complete}
    There is a finite number of 3-dicritical semi-complete digraphs.
\end{restatable}

While the proof of Theorem~\ref{thm:finite_3_dicritical_semi_complete} is fully human readable, the result is obtained by showing that the number of vertices of any $3$-dicritical semi-complete digraph does not exceed a pretty large number which originates from a Ramsey-type argument.

We after use a computer-assisted proof to provide the following characterization of all 3-dicritical semi-complete digraphs. 

\begin{restatable}{theorem}{thmallthreedicriticalsemicomplete}
\label{thm:all_3_dicritical_semi_complete}
    There are exactly eight 3-dicritical semi-complete digraphs. They are depicted in Figure~\ref{fig:all_3_dicritical}.
\end{restatable}

\newcolumntype{C}[1]{>{\centering\arraybackslash}p{#1}}

\begin{figure}[hbt!]
\begin{tabular}{C{.3\textwidth}C{.3\textwidth}C{.3\textwidth}}
%%%%%%%%%%%%%%%%%%%%%% K3 %%%%%%%%%%%%%%%%%%%%%%%%%%%%%%%%%%%
\subfigure [$\bid{K_3}$] {
    \begin{tikzpicture}[thick, every node/.style={transform shape}]
    	  \tikzset{vertex/.style = {circle,fill=black,minimum size=5pt, inner sep=0pt}}
            \tikzset{edge/.style = {->,> = latex'}}
            
        \draw[draw=white] (-1.5,-1) rectangle ++(3,2.8);
            \foreach \i in {0,...,2}{
                \node[vertex] (u\i) at (\i*360/3-30:1) {};
                \ifthenelse{\i>0}{
                    \pgfmathtruncatemacro{\j}{\i-1}
                    \foreach \k in {0,...,\j}{
                        \draw[edge, bend right=15] (u\i) to (u\k) {};
                        \draw[edge, bend right=15] (u\k) to (u\i) {};
                    }
                }{}
            }
    \end{tikzpicture}
} & 
%%%%%%%%%%%%%%%%%%%%%% W_3 %%%%%%%%%%%%%%%%%%%%%%%%%%%%%%%%%%%
\subfigure [$\ori{W_3}$] {
    \begin{tikzpicture}[thick, every node/.style={transform shape}]
    	  \tikzset{vertex/.style = {circle,fill=black,minimum size=5pt, inner sep=0pt}}
            \tikzset{edge/.style = {->,> = latex'}}

        \draw[draw=white] (-1.5,-1) rectangle ++(3,2.8);
            \node[vertex] (r) at (0,0) {};
            \foreach \i in {0,...,2}{
                \node[vertex] (u\i) at (\i*360/3-30:1) {};
                \draw[edge, bend right=15] (u\i) to (r) {};
                \draw[edge, bend right=15] (r) to (u\i) {};
            }
            \foreach \i in {0,...,2}{
                \pgfmathtruncatemacro{\j}{Mod(\i+1,3)}
                \draw[edge, bend right=30] (u\i) to (u\j) {};
            }
    \end{tikzpicture}
} &

%%%%%%%%%%%%%%%%%%%%%% K2 => K2 %%%%%%%%%%%%%%%%%%%%%%%%%%%%%%%%%%%
\subfigure [$\mathcal{R}(\bid{K_2},\bid{K_2})$] {
    \begin{tikzpicture}[thick,every node/.style={transform shape}]
        % Initialisation
        \tikzset{vertex/.style = {circle,fill=black,minimum size=5pt, inner sep=0pt}}
        \tikzset{edge/.style = {->,> = latex'}}

        \draw[draw=white] (-1.5,-0.3) rectangle ++(3,2.8);

        \node[vertex] (r) at (0,0) {};
        \node[] (arrow) at (0,1.6){\LARGE $\Ra$};
        \node[rotate=-135] (arrow) at (0.6,0.56){\LARGE $\Ra$};
        \node[rotate=135] (arrow) at (-0.6,0.56){\LARGE $\Ra$};
        \begin{scope}[xshift=-1.2cm, yshift=1.732cm]
            \foreach \i in {0,...,1}{
                \node[vertex] (x\i) at (\i*360/2-90:0.6) {};
            }
            \foreach \i in {0,...,1}{
                \pgfmathtruncatemacro{\j}{Mod(\i+1,2)}
                \draw[edge, bend right=40] (x\i) to (x\j) {};
            }
        \end{scope}
        \begin{scope}[xshift=1.2cm, yshift=1.732cm]
            \foreach \i in {0,...,1}{
                \node[vertex] (y\i) at (\i*360/2-90:0.6) {};
            }
            \foreach \i in {0,...,1}{
                \pgfmathtruncatemacro{\j}{Mod(\i+1,2)}
                \draw[edge, bend right=40] (y\i) to (y\j) {};
            }
        \end{scope}
    \end{tikzpicture}
} \\

%%%%%%%%%%%%%%%%%%%%%% H_5 %%%%%%%%%%%%%%%%%%%%%%%%%%%%%%%%%%%
\subfigure [$\mathcal{H}_5$] {
    \begin{tikzpicture}[thick,every node/.style={transform shape}]
        % Initialisation
        \tikzset{vertex/.style = {circle,fill=black,minimum size=5pt, inner sep=0pt}}
        \tikzset{edge/.style = {->,> = latex'}}
        
        \draw[draw=white] (-1.5,-1) rectangle ++(3,2.8);
        \foreach \i in {0,...,4}{
            \node[vertex] (u\i) at (\i*360/5+90:1) {};
        }
        \draw[edge, bend left=15] (u0) to (u4) {};
        \draw[edge, bend left=15] (u0) to (u1) {};
        \draw[edge, bend left=15] (u4) to (u0) {};
        \draw[edge, bend left=15] (u1) to (u0) {};
        \draw[edge] (u4) to (u1) {};
        \draw[edge] (u1) to (u3) {};
        \draw[edge] (u1) to (u2) {};
        \draw[edge] (u3) to (u4) {};
        \draw[edge] (u2) to (u4) {};
        \draw[edge] (u3) to (u0) {};
        \draw[edge] (u0) to (u2) {};
        \draw[edge] (u2) to (u3) {};
    \end{tikzpicture}
} &

%%%%%%%%%%%%%%%%%%%%%% K2 => C3 %%%%%%%%%%%%%%%%%%%%%%%%%%%%%%%%%%%
\subfigure [$\mathcal{R}(\bid{K_2},\ori{C_3})$] {
    \begin{tikzpicture}[thick,every node/.style={transform shape}]
        % Initialisation
        \tikzset{vertex/.style = {circle,fill=black,minimum size=5pt, inner sep=0pt}}
        \tikzset{edge/.style = {->,> = latex'}}

        \draw[draw=white] (-1.5,-0.3) rectangle ++(3,2.8);
        \node[vertex] (r) at (0,0) {};
        \node[] (arrow) at (0,1.6){\LARGE $\Ra$};
        \node[rotate=-135] (arrow) at (0.6,0.56){\LARGE $\Ra$};
        \node[rotate=135] (arrow) at (-0.6,0.56){\LARGE $\Ra$};
        \begin{scope}[xshift=-1.2cm, yshift=1.732cm]
            \foreach \i in {0,...,1}{
                \node[vertex] (x\i) at (\i*360/2-90:0.6) {};
            }
            \foreach \i in {0,...,1}{
                \pgfmathtruncatemacro{\j}{Mod(\i+1,2)}
                \draw[edge, bend right=40] (x\i) to (x\j) {};
            }
        \end{scope}
        \begin{scope}[xshift=1.2cm, yshift=1.732cm]
            \foreach \i in {0,...,2}{
                \node[vertex] (x\i) at (\i*360/3-90:0.6) {};
            }
            \foreach \i in {0,...,2}{
                \pgfmathtruncatemacro{\j}{Mod(\i+1,3)}
                \draw[edge, bend right=15] (x\i) to (x\j) {};
            }
        \end{scope}
    \end{tikzpicture}
} & 

%%%%%%%%%%%%%%%%%%%%%% C3 => K2 %%%%%%%%%%%%%%%%%%%%%%%%%%%%%%%%%%%
\subfigure [$\mathcal{R}(\ori{C_3},\bid{K_2})$] {
    \begin{tikzpicture}[thick,every node/.style={transform shape}]
        % Initialisation
        \tikzset{vertex/.style = {circle,fill=black,minimum size=5pt, inner sep=0pt}}
        \tikzset{edge/.style = {->,> = latex'}}

        \draw[draw=white] (-1.5,-0.3) rectangle ++(3,2.8);
        \node[vertex] (r) at (0,0) {};
        \node[] (arrow) at (0,1.6){\LARGE $\Ra$};
        \node[rotate=-135] (arrow) at (0.6,0.56){\LARGE $\Ra$};
        \node[rotate=135] (arrow) at (-0.6,0.56){\LARGE $\Ra$};
        \begin{scope}[xshift=-1.2cm, yshift=1.732cm]
            \foreach \i in {0,...,2}{
                \node[vertex] (x\i) at (\i*360/3-90:0.6) {};
            }
            \foreach \i in {0,...,2}{
                \pgfmathtruncatemacro{\j}{Mod(\i+1,3)}
                \draw[edge, bend right=15] (x\i) to (x\j) {};
            }
        \end{scope}
        \begin{scope}[xshift=1.2cm, yshift=1.732cm]
            \foreach \i in {0,...,1}{
                \node[vertex] (x\i) at (\i*360/2-90:0.6) {};
            }
            \foreach \i in {0,...,1}{
                \pgfmathtruncatemacro{\j}{Mod(\i+1,2)}
                \draw[edge, bend right=40] (x\i) to (x\j) {};
            }
        \end{scope}
    \end{tikzpicture}
} \\
\end{tabular}
\begin{tabular}{C{.1\textwidth}C{.3\textwidth}C{0.09\textwidth}C{.3\textwidth}}
&
%%%%%%%%%%%%%%%%%%%%%% C3 => C3 %%%%%%%%%%%%%%%%%%%%%%%%%%%%%%%%%%%
\subfigure [$\mathcal{R}(\ori{C_3},\ori{C_3})$] {
    \begin{tikzpicture}[thick,every node/.style={transform shape}]
        % Initialisation
        \tikzset{vertex/.style = {circle,fill=black,minimum size=5pt, inner sep=0pt}}
        \tikzset{edge/.style = {->,> = latex'}}

        \draw[draw=white] (-1.5,-0.3) rectangle ++(3,3.1);
        \node[vertex] (r) at (0,0) {};
        \node[] (arrow) at (0,1.6){\LARGE $\Ra$};
        \node[rotate=-135] (arrow) at (0.6,0.56){\LARGE $\Ra$};
        \node[rotate=135] (arrow) at (-0.6,0.56){\LARGE $\Ra$};
        \begin{scope}[xshift=-1.2cm, yshift=1.732cm]
            \foreach \i in {0,...,2}{
                \node[vertex] (x\i) at (\i*360/3-90:0.6) {};
            }
            \foreach \i in {0,...,2}{
                \pgfmathtruncatemacro{\j}{Mod(\i+1,3)}
                \draw[edge, bend right=15] (x\i) to (x\j) {};
            }
        \end{scope}
        \begin{scope}[xshift=1.2cm, yshift=1.732cm]
            \foreach \i in {0,...,2}{
                \node[vertex] (y\i) at (\i*360/3-90:0.6) {};
            }
            \foreach \i in {0,...,2}{
                \pgfmathtruncatemacro{\j}{Mod(\i+1,3)}
                \draw[edge, bend right=15] (y\i) to (y\j) {};
            }
        \end{scope}
    \end{tikzpicture}
} &  &
%%%%%%%%%%%%%%%%%%%%%% P_7 %%%%%%%%%%%%%%%%%%%%%%%%%%%%%%%%%%%
\subfigure [$\mathcal{P}_7$] {
    \begin{tikzpicture}[thick, every node/.style={transform shape}]
    	  \tikzset{vertex/.style = {circle,fill=black,minimum size=5pt, inner sep=0pt}}
            \tikzset{edge/.style = {->,> = latex'}}
            
            \draw[draw=white] (-1.5,-1.6) rectangle ++(3,3.1);

            \foreach \i in {0,...,6}{
                \node[vertex] (v\i) at (-51.428*\i-90:1.5) {};
            }
            \foreach \i in {0,...,6}{
                \pgfmathtruncatemacro{\j}{Mod(\i + 1,7)}
                \pgfmathtruncatemacro{\k}{Mod(\i + 2,7)}
                \pgfmathtruncatemacro{\l}{Mod(\i + 4,7)}
                \draw[edge] (v\i) to (v\j) {};
                \draw[edge] (v\i) to (v\k) {};
                \draw[edge] (v\i) to (v\l) {};
            }
    \end{tikzpicture}
}\\
\end{tabular}
\caption{The $3$-dicritical semi-complete digraphs, namely the bidirected complete graph $\bid{K_3}$, the directed wheel $\ori{W_3}$, the digraph $\mathcal{H}_5$, the rotative digraphs $\mathcal{R}(H_1,H_2)$ for every $H_1,H_2\in \{\bid{K_2},\ori{C_3}\}$ and the Paley tournament on seven vertices $\mathcal{P}_7$. A big arrow linking two sets of vertices indicates that there is exactly one arc from every vertex in the first set to every vertex in the second set.}
\label{fig:all_3_dicritical}
\end{figure}

In particular, we can characterize all 3-dicritical tournaments. 

\begin{corollary}
    There are exactly two 3-dicritical tournaments, namely $\mathcal{R}(\ori{C_3},\ori{C_3})$ and $\mathcal{P}_7$.
\end{corollary}

We finally investigate the maximum density of $3$-dicritical digraphs. The {\it bidirected part} of a digraph $D$ is the graph $B(D)$ with vertex set $V(D)$ in which two vertices are linked by an edge if and only if there is a digon between them in $D$. We prove in Proposition~\ref{prop:digon_forest} that $B(D)$ is a forest for every $3$-dicritical digraph $D$ that is not a bidirected odd cycle. From this result, one can easily deduce that $m(D) \leq \binom{n}{2} + n-1$ holds for every $3$-dicritical digraph $D$ different from $\bid{K_3}$. We slightly improve this upper bound on $m(D)$ as follows (the digraph $\ori{W_3}$ is depicted in Figure~\ref{fig:all_3_dicritical}).

\begin{restatable}{theorem}{thmboundtwothird}
    \label{thm:bound_two_third}
    If $D$ is a $3$-dicritical digraph distinct from $\bid{K_3}$ and $\ori{W_3}$, then $m(D)\leq  \binom{n}{2} + \frac{2}{3}n$. 
\end{restatable}

The rest of this article is structured as follows: in Section~\ref{prel}, we give a collection of preliminary results which will be used in the proofs of Theorems~\ref{thm:finite_3_dicritical_semi_complete},~\ref{thm:all_3_dicritical_semi_complete}, and~\ref{thm:bound_two_third}. In Section~\ref{secfaible}, we prove Theorem~\ref{thm:finite_3_dicritical_semi_complete}. In Section~\ref{secfort}, we prove Theorem~\ref{thm:all_3_dicritical_semi_complete}, with all code we use being shifted to the appendix. Section~\ref{sec:bound_two_third} is devoted to the proof of Theorem~\ref{thm:bound_two_third}.
Finally, in Section~\ref{conc}, we conclude our work and give some directions for further research. 

\section{Useful lemmas}\label{prel}

In this section, we give a collection of preliminary results we need in the proof of Theorem~\ref{thm:finite_3_dicritical_semi_complete}. Most of them will be reused in the proof of Theorem~\ref{thm:all_3_dicritical_semi_complete}. We first describe 2-dicolourings with some important extra properties.
\medskip

Let $D$ be a digraph and $uv$ be an arc of $D$. A \textit{$uv$-colouring} of $D$ is a $2$-colouring $\phi : V(D) \xrightarrow[]{} [2]$ such that:
\begin{itemize}
    \item $\phi$ is a $2$-dicolouring of $D\setminus uv$,
    \item $\phi(u) = \phi(v) = 1$, and
    \item $D\setminus uv$, coloured with $\phi$, does not contain any monochromatic directed $uv$-path.
\end{itemize}

There is a close relationship between $3$-dicritical digraphs and $uv$-colourings.
\begin{lemma}\label{lemma:dicol}
    Let $D$ be a 3-dicritical digraph and $uv$ be an arc of $D$. Then $D$ admits a $uv$-colouring.
\end{lemma}
\begin{proof}
    As $D$ is 3-dicritical, there is a 2-dicolouring $\phi:V(D)\rightarrow [2]$ of $D\setminus uv$. By symmetry, we may suppose $\phi(u)=1$. As $\phi$ is not a 2-dicolouring of $D$, we obtain that $\phi(v)=1$ and there is a directed $vu$-path $P$ in $D$ such that $\phi(x)=1$ for all $x \in V(P)$. If there is also a directed $uv$-path $Q$ in $D\setminus uv$ such that $\phi(x)=1$ for all $x \in V(Q)$, then the subdigraph of $D\setminus uv$ induced by $V(P)\cup V(Q)$ contains a monochromatic directed cycle. This contradicts $\phi$ being a 2-dicolouring of $D\setminus uv$.
\end{proof}

The next result showing that every arc of a 3-dicritical semi-complete digraph is contained in a short directed cycle will play a crucial role in the upcoming proofs.
\begin{lemma}\label{lemma:every_arc_in_triangle}
    Let $D$ be a 3-dicritical semi-complete digraph. Then every arc $a \in A(D)$ either belongs to a digon or is contained in an induced directed triangle.
\end{lemma}
\begin{proof}
     As $D$ is 3-dicritical, there is a  $2$-dicolouring $\phi$ of $D\setminus a$. As $\phi$ is not a 2-dicolouring of $D$, there exists a directed cycle $C$ in $D$ such that $C$ is monochromatic with respect to $\phi$. We may suppose that $C$ is chosen to be of minimum length with this property. As $D$ is semi-complete, we obtain that $C$ is either a digon or an induced directed triangle. As $C$ is not a monochromatic directed cycle of $D\setminus a$ with respect to $\phi$, we obtain that $a \in A(C)$.
\end{proof}

We define $O_5$ as the oriented graph which consists of a directed triangle $xyz$ and two additional vertices $u,v$, one arc from $u$ to every vertex of the directed triangle, one arc from every vertex of the directed triangle to $v$, and the arc $uv$. An illustration can be found in Figure~\ref{fig:o5}.

\begin{figure}[hbt!]
    \begin{center}	
          \begin{tikzpicture}[thick,scale=1, every node/.style={transform shape}]
            \tikzset{vertex/.style = {circle,fill=black,minimum size=5pt, inner sep=0pt}}
            \tikzset{edge/.style = {->,> = latex'}}

            \node[vertex, label=above:$x$] (x) at (30:0.7) {};
            \node[vertex, label=above:$y$] (y) at (150:0.7) {};
            \node[vertex, label=below:$z$] (z) at (-90:0.7) {}; 
            \node[vertex, label=below:$u$] (u) at (-135:2.1) {};
            \node[vertex, label=below:$v$] (v) at (-45:2.1) {}; 

            \draw[edge] (u) to (v) {};
            \draw[edge] (u) to (x) {};
            \draw[edge] (u) to (y) {};
            \draw[edge] (u) to (z) {};
            \draw[edge] (x) to (v) {};
            \draw[edge] (y) to (v) {};
            \draw[edge] (z) to (v) {};
            \draw[edge] (x) to (y) {};
            \draw[edge] (y) to (z) {};
            \draw[edge] (z) to (x) {};
          \end{tikzpicture}
      \caption{The oriented graph $O_5$.}
      \label{fig:o5}
    \end{center}
\end{figure}

The following result is a consequence of Lemma~\ref{lemma:dicol}.
\begin{lemma}\label{lemma:o5_free}
    Let $D$ be a 3-dicritical digraph. Then $D$ does not contain $O_5$ as a subdigraph.
\end{lemma}
\begin{proof}
    Assume for a contradiction that $D$ contains $O_5$ as a subdigraph and let $V(O_5)=\{u,v,x,y,z\}$ be the labelling depicted in Figure~\ref{fig:o5}. By Lemma~\ref{lemma:dicol}, there exists a $uv$-colouring $\phi$ of $D$. Since there exists no monochromatic directed $uv$-path, we have $\phi(x)=\phi(y)=\phi(z)=2$. Hence $D\setminus uv$ contains a monochromatic directed triangle with respect to $\phi$, a contradiction.
\end{proof}

Let $S$ be a transitive subtournament of a digraph $D=(V,A)$. We denote by $v_1,\dots,v_s$ the unique acyclic ordering of $S$.
For some $i,j \in [s]$, we say that $\{v_i,\ldots,v_j\}$ is an {\it interval} of $S$. Observe that $\emptyset$ is an interval. 
For $i_0,j_0,i_1,j_1\in [s]$ with $j_0<i_1$, we say that the interval $\{v_{i_0},\ldots,v_{j_0}\}$ is {\it smaller} than the interval $\{v_{i_1},\ldots,v_{j_1}\}$. By convention $\emptyset$ is both smaller and greater than any other interval. A sequence of intervals $P_1,\ldots,P_t$ is called {\it increasing} if $P_i$ is smaller than $P_j$ for all $i,j \in [t]$ with $i<j$.

\begin{lemma}\label{lemma:intervals}
    Let $T$ be a subtournament of a $3$-dicritical digraph $D$ and let $S$ be a transitive subtournament of $T$ with acyclic ordering $v_1,\dots,v_s$.
     For any $x \in V(T) \setminus S$, there is an increasing sequence of intervals $(I_1, I_2, I_3, I_4)$ with $\bigcup_{i=1}^4 I_i=S$ such that, in $T$, $x$ dominates $I_1 \cup I_3$ and is dominated by $I_2 \cup I_4$ .      
\end{lemma}
 \begin{proof}
    Assume this is not the case.
    Then there exists an increasing sequence of indices $(i_1 , i_2,  i_3,  i_4)$ such that, in $T$, $x$ is dominated by $v_{i_1}$ and $v_{i_3}$ and dominates $v_{i_2}$ and $v_{i_4}$.
    Then the subdigraph of $D$ induced by $\{v_{i_1}, v_{i_2},x, v_{i_3}, v_{i_4}\}$ contains $O_5$ as a subdigraph, a contradiction to Lemma~\ref{lemma:o5_free}.
 \end{proof}

 We finally need a well-known theorem which can be found in many basic textbooks on graph theory, see for example~\cite[Theorem~9.1.3]{diestel2017}.

 \begin{theorem}[{\sc Multi-colour Ramsey Theorem}]\label{ramsey}
 Let $a$ and $b$ be positive integers. There exists a smallest integer $R_a(b)$ such that for $G$ being a copy of $K_{R_a(b)}$ and for every mapping $\psi:E(G)\rightarrow [a]$, there is a set $S \subseteq V(G)$ of cardinality $b$ and $i \in [a]$ such that $\psi(e)=i$ for all $e \in E(G\ind{S})$.
 \end{theorem}

\section{A simple proof for finiteness}\label{secfaible}

In this section, we prove that the number of $3$-dicritical semi-complete digraphs is finite. Let us first restate this result.

\thmfinitethreedicriticalsemicomplete*

\begin{proof}
    Let $D=(V,A)$ be a $3$-dicritical semi-complete digraph. We will show that $n(D) \leq 12R_6(3) + 1$, where $R_6(3)$ refers to the Ramsey number in Theorem~\ref{ramsey}. Assume for the sake of a contradiction that  $n(D) \geq 12R_6(3) + 2$. 
    Let $S\subseteq V$ be a maximum set of vertices such that $D\ind{S}$ is acyclic. Let $v_1,\dots,v_s$ be the unique acyclic ordering of $S$. 
    Since $D$ is $3$-dicritical, for an arbitrary vertex $x \in V$, we have that $D- x$ is 2-dicolourable. This yields $s\geq \left\lceil \frac{n(D)-1}{2} \right\rceil \geq 6R_6(3)+1$. 
    
    By Lemma~\ref{lemma:every_arc_in_triangle}, for every $i\in[s-1]$, the arc $v_iv_{i+1}$ belongs to a digon or an induced directed triangle. Therefore, since $D\ind{S}$ is acyclic, we know that there exists a vertex $x_i\in V\setminus S$ such that $v_iv_{i+1}x_iv_i$ is an induced  directed triangle $C_i$.

    Let $T$ be an arbitrary spanning subtournament of $D$. Observe that $T\ind{S}=D\ind{S}$ as $D\ind{S}$ is acyclic. Further, the directed triangle $C_i$ is contained in $T$ for $i \in [s-1]$ as $C_i$ is induced in $D$.
    For any vertex $x$ in $V\setminus S$ and $i\in [s-1]$, we say that $x$ {\it switches at $i$} if $x$ dominates $v_i$ and is dominated by $v_{i+1}$ in $T$ or $x$ is dominated by $v_i$ and dominates $v_{i+1}$ in $T$.
    
    Let $H$ be the digraph with vertex set $V(H)=[s-1]$ and arc set $A(H) = A_1\cup A_2$ with $A_1=\{(i,i+1) \mid i\in [s-2]\}$ and $A_2=\{(i,j) \mid i\neq j \mbox{ and $x_i$ switches at $j$.}\}$.

   By Lemma~\ref{lemma:intervals}, for $i \in [s-1]$, we have that $x_i$ switches at 
at most three indices in $[s-1]$. Further, as $C_i$ is a directed triangle, $x_i$ switches at $i$ which yields that $x_i$ switches at at most two indices in $[s-1]\setminus \{i\}$.
Thus every $i \in [s-1]$ is the tail of at most two arcs in $A_2$.

    For every subset $J$ of $[s-1]$, observe that $H\ind{J}$ contains  at most $|J|-1$ arcs in $A_1$ and at most $2|J|$ arcs in $A_2$, hence at most $3|J| -1$ arcs in total. Thus $\UG(H)\ind{J}$ has a vertex of degree at most $5$.
    Hence $\UG(H)$ is $5$-degenerate, and so it is $6$-colourable. Therefore $H$ has an independent set $I$ of size $\left\lceil \frac{1}{6}(s-1) \right\rceil \geq R_6(3)$.

      By definition of $I$, for any $i,j \in I$ with $i\neq j$, we have that  $V(C_i)$ and $V(C_j)$ are disjoint. Moreover, either $\{v_j,v_{j+1}\}$ dominates $x_i$ in $T$ or  $\{v_j,v_{j+1}\}$ is dominated by $x_i$ in $T$. Hence, if $i<j$, the subdigraph of $T$ induced by $V(C_i)\cup V(C_j)$ is  one of the eight tournaments depicted in Figure~\ref{fig:configurations_Ci_Cj}. For $(\alpha)\in \{(a),\ldots,(h)\}$, we say that $(i,j)$ is an {\it $(\alpha)$-configuration} if $T\ind{V(C_i)\cup V(C_j)}$ is the tournament depicted in Figure~\ref{fig:configurations_Ci_Cj}~$(\alpha)$. 

    \begin{figure}[hbt!]
    \begin{center}	
      \begin{tikzpicture}[thick,scale=1, every node/.style={transform shape}]
        \tikzset{vertex/.style = {circle,fill=black,minimum size=5pt, inner sep=0pt}}
        \tikzset{edge/.style = {->,> = latex'}}
        \begin{scope}
            \node[vertex, label=left:$v_{i+1}$] (vip1) at (0:0.45) {};
            \node[vertex, label=left:$v_{i}$] (vi) at (120:1) {};
            \node[vertex, label=left:$x_{i}$] (xi) at (-120:1) {};
            \draw[edge] (xi) to (vi);
            \draw[edge] (vi) to (vip1);
            \draw[edge] (vip1) to (xi);
            
            \begin{scope}[xshift=3cm]
            \node[vertex, label=right:$v_{j+1}$] (vjp1) at (0:0.45) {};
            \node[vertex, label=right:$v_{j}$] (vj) at (120:1) {};
            \node[vertex, label=right:$x_{j}$] (xj) at (-120:1) {};
            \draw[edge] (xj) to (vj);
            \draw[edge] (vj) to (vjp1);
            \draw[edge] (vjp1) to (xj);
            \end{scope}

            \draw[edge,dashed] (xi) to (xj);
            \draw[edge,dashed] (vip1) to (vj);
            \draw[edge,dashed] (vip1) to (vjp1);
            \draw[edge,dashed] (vi) to (vj);
            \draw[edge,dashed] (vi) to (vjp1);
            
            \draw[edge] (vj) to (xi);
            \draw[edge] (vjp1) to (xi);
            \draw[edge] (xj) to (vi);
            \draw[edge] (xj) to (vip1);

            \node[] (a) at (1.3,-1.3) {$(a)$};
        \end{scope}
        \begin{scope}[xshift=7cm]
            \node[vertex, label=left:$v_{i+1}$] (vip1) at (0:0.45) {};
            \node[vertex, label=left:$v_{i}$] (vi) at (120:1) {};
            \node[vertex, label=left:$x_{i}$] (xi) at (-120:1) {};
            \draw[edge] (xi) to (vi);
            \draw[edge] (vi) to (vip1);
            \draw[edge] (vip1) to (xi);
            
            \begin{scope}[xshift=3cm]
            \node[vertex, label=right:$v_{j+1}$] (vjp1) at (0:0.45) {};
            \node[vertex, label=right:$v_{j}$] (vj) at (120:1) {};
            \node[vertex, label=right:$x_{j}$] (xj) at (-120:1) {};
            \draw[edge] (xj) to (vj);
            \draw[edge] (vj) to (vjp1);
            \draw[edge] (vjp1) to (xj);
            \end{scope}

            \draw[edge,dashed] (vip1) to (vj);
            \draw[edge,dashed] (vip1) to (vjp1);
            \draw[edge,dashed] (vi) to (vj);
            \draw[edge,dashed] (vi) to (vjp1);
            
            \draw[edge] (xj) to (xi);
            \draw[edge] (vj) to (xi);
            \draw[edge] (vjp1) to (xi);
            \draw[edge] (xj) to (vi);
            \draw[edge] (xj) to (vip1);

            \node[] (a) at (1.3,-1.3) {$(b)$};
        \end{scope}
        \begin{scope}[yshift=-3cm]
            \node[vertex, label=left:$v_{i+1}$] (vip1) at (0:0.45) {};
            \node[vertex, label=left:$v_{i}$] (vi) at (120:1) {};
            \node[vertex, label=left:$x_{i}$] (xi) at (-120:1) {};
            \draw[edge] (xi) to (vi);
            \draw[edge] (vi) to (vip1);
            \draw[edge] (vip1) to (xi);
            
            \begin{scope}[xshift=3cm]
            \node[vertex, label=right:$v_{j+1}$] (vjp1) at (0:0.45) {};
            \node[vertex, label=right:$v_{j}$] (vj) at (120:1) {};
            \node[vertex, label=right:$x_{j}$] (xj) at (-120:1) {};
            \draw[edge] (xj) to (vj);
            \draw[edge] (vj) to (vjp1);
            \draw[edge] (vjp1) to (xj);
            \end{scope}

            \draw[edge,dashed] (vip1) to (vj);
            \draw[edge,dashed] (vip1) to (vjp1);
            \draw[edge,dashed] (vi) to (vj);
            \draw[edge,dashed] (vi) to (vjp1);
            \draw[edge,dashed] (xi) to (vj);
            \draw[edge,dashed] (xi) to (xj);
            \draw[edge,dashed] (xi) to (vjp1);
            
            \draw[edge] (xj) to (vi);
            \draw[edge] (xj) to (vip1);

            \node[] (a) at (1.3,-1.3) {$(c)$};
        \end{scope}
        \begin{scope}[yshift=-3cm, xshift=7cm]
            \node[vertex, label=left:$v_{i+1}$] (vip1) at (0:0.45) {};
            \node[vertex, label=left:$v_{i}$] (vi) at (120:1) {};
            \node[vertex, label=left:$x_{i}$] (xi) at (-120:1) {};
            \draw[edge] (xi) to (vi);
            \draw[edge] (vi) to (vip1);
            \draw[edge] (vip1) to (xi);
            
            \begin{scope}[xshift=3cm]
            \node[vertex, label=right:$v_{j+1}$] (vjp1) at (0:0.45) {};
            \node[vertex, label=right:$v_{j}$] (vj) at (120:1) {};
            \node[vertex, label=right:$x_{j}$] (xj) at (-120:1) {};
            \draw[edge] (xj) to (vj);
            \draw[edge] (vj) to (vjp1);
            \draw[edge] (vjp1) to (xj);
            \end{scope}

            \draw[edge,dashed] (xi) to (xj);
            \draw[edge,dashed] (vip1) to (vj);
            \draw[edge,dashed] (vip1) to (vjp1);
            \draw[edge,dashed] (vi) to (vj);
            \draw[edge,dashed] (vi) to (vjp1);
            \draw[edge,dashed] (vi) to (xj);
            \draw[edge,dashed] (vip1) to (xj);
            
            \draw[edge] (vj) to (xi);
            \draw[edge] (vjp1) to (xi);

            \node[] (a) at (1.3,-1.3) {$(d)$};
        \end{scope}
        \begin{scope}[yshift=-6cm]
            \node[vertex, label=left:$v_{i+1}$] (vip1) at (0:0.45) {};
            \node[vertex, label=left:$v_{i}$] (vi) at (120:1) {};
            \node[vertex, label=left:$x_{i}$] (xi) at (-120:1) {};
            \draw[edge] (xi) to (vi);
            \draw[edge] (vi) to (vip1);
            \draw[edge] (vip1) to (xi);
            
            \begin{scope}[xshift=3cm]
            \node[vertex, label=right:$v_{j+1}$] (vjp1) at (0:0.45) {};
            \node[vertex, label=right:$v_{j}$] (vj) at (120:1) {};
            \node[vertex, label=right:$x_{j}$] (xj) at (-120:1) {};
            \draw[edge] (xj) to (vj);
            \draw[edge] (vj) to (vjp1);
            \draw[edge] (vjp1) to (xj);
            \end{scope}

            \draw[edge,dashed] (vip1) to (vj);
            \draw[edge,dashed] (vip1) to (vjp1);
            \draw[edge,dashed] (vi) to (vj);
            \draw[edge,dashed] (vi) to (vjp1);
            \draw[edge,dashed] (vi) to (xj);
            \draw[edge,dashed] (vip1) to (xj);
            \draw[edge,dashed] (xi) to (vj);
            \draw[edge,dashed] (xi) to (vjp1);

            \draw[edge] (xj) to (xi);
            
            \node[] (a) at (1.3,-1.3) {$(e)$};
        \end{scope}
        \begin{scope}[xshift=7cm,yshift=-6cm]
            \node[vertex, label=left:$v_{i+1}$] (vip1) at (0:0.45) {};
            \node[vertex, label=left:$v_{i}$] (vi) at (120:1) {};
            \node[vertex, label=left:$x_{i}$] (xi) at (-120:1) {};
            \draw[edge] (xi) to (vi);
            \draw[edge] (vi) to (vip1);
            \draw[edge] (vip1) to (xi);
            
            \begin{scope}[xshift=3cm]
            \node[vertex, label=right:$v_{j+1}$] (vjp1) at (0:0.45) {};
            \node[vertex, label=right:$v_{j}$] (vj) at (120:1) {};
            \node[vertex, label=right:$x_{j}$] (xj) at (-120:1) {};
            \draw[edge] (xj) to (vj);
            \draw[edge] (vj) to (vjp1);
            \draw[edge] (vjp1) to (xj);
            \end{scope}

            \draw[edge,dashed] (vip1) to (vj);
            \draw[edge,dashed] (vip1) to (vjp1);
            \draw[edge,dashed] (vi) to (vj);
            \draw[edge,dashed] (vi) to (vjp1);
            \draw[edge,dashed] (vi) to (xj);
            \draw[edge,dashed] (vip1) to (xj);
            \draw[edge,dashed] (xi) to (vj);
            \draw[edge,dashed] (xi) to (vjp1);
            \draw[edge,dashed] (xi) to (xj);
            \node[] (a) at (1.3,-1.3) {$(f)$};
        \end{scope}
        \begin{scope}[yshift=-9cm]
            \node[vertex, label=left:$v_{i+1}$] (vip1) at (0:0.45) {};
            \node[vertex, label=left:$v_{i}$] (vi) at (120:1) {};
            \node[vertex, label=left:$x_{i}$] (xi) at (-120:1) {};
            \draw[edge] (xi) to (vi);
            \draw[edge] (vi) to (vip1);
            \draw[edge] (vip1) to (xi);
            
            \begin{scope}[xshift=3cm]
            \node[vertex, label=right:$v_{j+1}$] (vjp1) at (0:0.45) {};
            \node[vertex, label=right:$v_{j}$] (vj) at (120:1) {};
            \node[vertex, label=right:$x_{j}$] (xj) at (-120:1) {};
            \draw[edge] (xj) to (vj);
            \draw[edge] (vj) to (vjp1);
            \draw[edge] (vjp1) to (xj);
            \end{scope}

            \draw[edge,dashed] (vip1) to (vj);
            \draw[edge,dashed] (vip1) to (vjp1);
            \draw[edge,dashed] (vi) to (vj);
            \draw[edge,dashed] (vi) to (vjp1);
            \draw[edge,dashed] (vi) to (xj);
            \draw[edge,dashed] (vip1) to (xj);
            
            \draw[edge] (xj) to (xi);
            \draw[edge] (vj) to (xi);
            \draw[edge] (vjp1) to (xi);
            
            \node[] (a) at (1.3,-1.3) {$(g)$};
        \end{scope}
        \begin{scope}[xshift=7cm,yshift=-9cm]
        \node[vertex, label=left:$v_{i+1}$] (vip1) at (0:0.45) {};
            \node[vertex, label=left:$v_{i}$] (vi) at (120:1) {};
            \node[vertex, label=left:$x_{i}$] (xi) at (-120:1) {};
            \draw[edge] (xi) to (vi);
            \draw[edge] (vi) to (vip1);
            \draw[edge] (vip1) to (xi);
            
            \begin{scope}[xshift=3cm]
            \node[vertex, label=right:$v_{j+1}$] (vjp1) at (0:0.45) {};
            \node[vertex, label=right:$v_{j}$] (vj) at (120:1) {};
            \node[vertex, label=right:$x_{j}$] (xj) at (-120:1) {};
            \draw[edge] (xj) to (vj);
            \draw[edge] (vj) to (vjp1);
            \draw[edge] (vjp1) to (xj);
            \end{scope}

            \draw[edge,dashed] (vip1) to (vj);
            \draw[edge,dashed] (vip1) to (vjp1);
            \draw[edge,dashed] (vi) to (vj);
            \draw[edge,dashed] (vi) to (vjp1);
            \draw[edge,dashed] (xi) to (vj);
            \draw[edge,dashed] (xi) to (vjp1);
            
            \draw[edge] (xj) to (vi);
            \draw[edge] (xj) to (vip1);
            \draw[edge] (xj) to (xi);
            
            \node[] (a) at (1.3,-1.3) {$(h)$};
        \end{scope}
      \end{tikzpicture}
      \caption{A listing of all possible configurations for $i,j \in I$ with $i<j$. For the sake of better readability, the arcs in $A(C_i)\cup A(C_j)$ and the arcs from $V(C_j)$ to $V(C_i)$ are solid, and the arcs from $V(C_i)$ to $V(C_j)$ are dashed.}\label{fig:configurations_Ci_Cj}
    \end{center}
    \end{figure}

    Let us fix a pair $i,j\in I$ with $i<j$. 
    We know that it is not a $(g)$-configuration, for otherwise $D\ind{v_{i+1},v_j,v_{j+1},x_j,x_i}$ contains  $O_5$ as a subdigraph, a contradiction to Lemma~\ref{lemma:o5_free}. We also know that it is not an $(h)$-configuration, for otherwise $D\ind{x_j,v_i,v_{i+1},x_i,v_j}$ contains $O_5$ as a subdigraph, a contradiction to Lemma~\ref{lemma:o5_free}.

    Since $|I|\geq R_6(3)$, and by definition of $R_6(3)$, we know that there exist $\{ i,j,h \} \subseteq I$, $i<j<h$, and $(\alpha) \in \{(a),\dots,(e)\}$ such that the three pairs $(i,j),(j,h),(i,h)$ are $(\alpha)$-configurations. We show that each of the six cases yields a contradiction, implying the result.
    \begin{itemize}
        \item If $(\alpha) = (a)$, let $\phi$ be a $v_{i+1}v_{h}$-colouring of $D$, the existence of which is guaranteed by Lemma~\ref{lemma:dicol}. Recall that $\phi(v_{i+1}) = \phi(v_{h}) = 1$, $\phi$ is a $2$-dicolouring of $D \setminus v_{i+1}v_{h}$ and $D$ coloured with $\phi$ contains no monochromatic directed $v_{i+1}v_{h}$-path. Then $\phi(v_j) = \phi(v_{j+1}) =  2$ because $\{v_j,v_{j+1}\} \subseteq N_D^+(v_{i+1}) \cap N_D^-(v_{h})$. Thus, since $C_j$ is not monochromatic in $\phi$, we have $\phi(x_j) = 1$. We obtain that $\phi(x_i) = \phi(v_i) = 2$, for otherwise $v_{i+1}x_ix_jv_{i+1}$ or $v_iv_hx_jv_i$ is monochromatic.
        We deduce that $v_iv_jx_iv_i$ is monochromatic, a contradiction. 

        \item If $(\alpha) = (b)$, then $D\ind{x_h,v_j,v_{j+1},x_j,x_i}$ contains $O_5$ as a subdigraph, a contradiction to Lemma~\ref{lemma:o5_free}.

        \item If $(\alpha) = (c)$, then $D\ind{x_i,v_j,v_{j+1},x_j,v_h}$ contains $O_5$ as a subdigraph, a contradiction to Lemma~\ref{lemma:o5_free}.
        
        \item If $(\alpha) = (d)$, then $D\ind{v_i,v_j,v_{j+1},x_j,x_h}$ contains $O_5$ as a subdigraph, a contradiction to Lemma~\ref{lemma:o5_free}.
        
        \item If $(\alpha) = (e)$, then $D\ind{v_i,v_j,v_{j+1},x_j,v_h}$ contains $O_5$ as a subdigraph, a contradiction to Lemma~\ref{lemma:o5_free}.

        \item If $(\alpha) = (f)$, then $D\ind{v_i,v_j,v_{j+1},x_j,v_h}$ contains $O_5$ as a subdigraph, a contradiction to Lemma~\ref{lemma:o5_free}.
    \end{itemize}
\end{proof}

\section{The 3-dicritical semi-complete digraphs}
\label{secfort}

This section is devoted to a computer-assisted proof of Theorem~\ref{thm:all_3_dicritical_semi_complete}. 
 It follows a similar line as the one of Theorem~\ref{thm:finite_3_dicritical_semi_complete}, but it needs some refined arguments. Further, due to the significant number of necessary computations, several parts of the proof are computer-assisted. We used codes implemented using SageMath. They are accessible on the third author's \href{https://github.com/lucaspicasarri/three_dicritical_semicomplete_digraphs}{GitHub page} and are given in the appendix.
 
 We first restrain the structure of 3-dicritical semi-complete digraphs. To prove Theorem~\ref{thm:finite_3_dicritical_semi_complete}, we only needed the fact that $O_5$ does not occur as a subdigraph. 
 To prove Theorem~\ref{thm:all_3_dicritical_semi_complete}, we need to prove that several other digraphs cannot be subdigraphs or induced subdigraphs of a 3-dicritical digraph. One of these digraphs is the transitive tournament of size at least 8.  
 While already parts of this proof are computer-assisted, the most intense computation part is carried out after. We generate all semi-complete digraphs satisfying these properties and check that none of them has dichromatic number 3, except the ones depicted in Figure~\ref{fig:all_3_dicritical}.

 Before dealing with the collection of digraphs which are not contained in 3-dicritical semi-complete digraphs as subdigraphs, we first give the following simple observation on matchings in graphs on seven vertices which will prove useful later on.

 \begin{lemma}\label{lemma:matchpath}
    Let $H$ be a  graph that is obtained from a path $w_1\ldots w_7$ by adding the edges of a matching $M$ on $\{w_1,\ldots,w_7\}$. Then there is a stable set $S\subseteq V(H)$ with $|S|=3$ and $\{w_1,w_7\}\setminus S\neq \emptyset$.
\end{lemma}
\begin{proof}
    Suppose otherwise. If $w_1$ is not incident to an edge of $M$, then, as none of $\{w_1,w_3,w_5\}$ and $\{w_1,w_3,w_6\}$ is an independent set, we obtain that $w_3w_5,w_3w_6 \in E(M)$, a contradiction to $M$ being a matching. Hence $M$ contains an edge $e_1$ incident to $w_1$. Similarly, $M$ contains an edge $e_7$ incident to $w_7$. Further, $M$ contains an edge $e_0$ both of whose endvertices are contained in $\{w_2,w_4,w_6\}$. As none of $e_0$ and $e_7$ are contained in one of $\{w_1,w_3,w_6\}$ and $\{w_1,w_3,w_5\}$, we obtain that $e_1$ needs to be contained in both of them. This yields $e_1=w_1w_3$. Similarly, we obtain $e_7=w_5w_7$. By symmetry, we may suppose that $e_0\neq w_4w_6$. But then $\{w_1,w_4,w_6\}$ is an independent set, a contradiction.
\end{proof}

We now start excluding some subdigraphs of $3$-dicritical semi-complete digraphs. 
We first define $O_4$ as the digraph which consists of a copy of $\bid{K_2}$ and two additional vertices $u,v$, one arc from $u$ to every vertex of $\bid{K_2}$, one arc from every vertex of $\bid{K_2}$ to $v$, and the arc $uv$. An illustration can be found in Figure~\ref{fig:o4}.

\begin{figure}[hbt!]
    \begin{center}	
          \begin{tikzpicture}[thick,scale=1, every node/.style={transform shape}]
            \tikzset{vertex/.style = {circle,fill=black,minimum size=5pt, inner sep=0pt}}
            \tikzset{edge/.style = {->,> = latex'}}

            \node[vertex, label=above:$x$] (x) at (90:0.7) {};
            \node[vertex, label=below:$y$] (z) at (-90:0.7) {}; 
            \node[vertex, label=below:$u$] (u) at (180:1.4) {};
            \node[vertex, label=below:$v$] (v) at (0:1.4) {}; 

            \draw[edge] (u) to (v) {};
            \draw[edge] (u) to (x) {};
            \draw[edge] (u) to (z) {};
            \draw[edge] (x) to (v) {};
            \draw[edge] (z) to (v) {};
            \draw[edge,bend left=15] (x) to (z) {};
            \draw[edge,bend left=15] (z) to (x) {};
          \end{tikzpicture}
      \caption{The digraph $O_4$.}
      \label{fig:o4}
    \end{center}
\end{figure}
The digraph $O_4$ plays a similar role as $O_5$. Also, the proof of the following result is similar to the one of Lemma~\ref{lemma:o5_free}.
\begin{lemma}\label{lemma:o4_free}
    Let $D$ be a 3-dicritical digraph. Then $D$ does not contain $O_4$ as a subdigraph.
\end{lemma}
\begin{proof}
    Assume for the purpose of contradiction that $D$ contains $O_4$ as a subdigraph and let $V(O_4)=\{u,v,x,y\}$ be the labelling depicted in Figure~\ref{fig:o4}. By Lemma~\ref{lemma:dicol}, there exists a 2-dicolouring $\phi$ of $D\setminus uv$ with $\phi(x)=\phi(y)$. Hence $D\setminus uv$ contains a monochromatic digon with respect to $\phi$, a contradiction.
\end{proof}
In the following, let $\bid{S_4}$ be the bidirected star on 4 vertices, see Figure~\ref{4etoile}. The following result shows that $\bid{S_4}$ cannot be the subdigraph of any large 3-dicritical semi-complete digraph.
 \begin{figure}[ht]\begin{center}
       \begin{tikzpicture}[thick, every node/.style={transform shape}]
    	  \tikzset{vertex/.style = {circle,fill=black,minimum size=5pt, inner sep=0pt}}
            \tikzset{edge/.style = {->,> = latex'}}

            \node[vertex] (r) at (0,0) {};
            \foreach \i in {0,...,2}{
                \node[vertex] (u\i) at (\i*360/3-30:1) {};
                \draw[edge, bend right=15] (u\i) to (r) {};
                \draw[edge, bend right=15] (r) to (u\i) {};
            }
    \end{tikzpicture}
      \caption{The bidirected star on 4 vertices $\bid{S_4}$.}\label{4etoile}
    \end{center}
 \end{figure}
\begin{lemma}
    \label{lemma:S4_free}
    Let $D$ be a semi-complete digraph containing $\bid{S_4}$ as a subdigraph. Then $D$ is 3-dicritical if and only if $D$ is $\ori{W_3}$.
\end{lemma}
\begin{proof}
    It is easy to see that $\ori{W_3}$ is 3-dicritical and contains $\bid{S_4}$. For the other direction, let $D$ be a $3$-dicritical semi-complete digraph such that $D$ contains a vertex $u$ linked by digons to three distinct vertices $x,y,z$.

    Then, as $D$ is semi-complete, we have that $D\ind{\{x,y,z\}}$ needs to contain $\ori{C_3}$ or $TT_3$ as a subdigraph. If it is $TT_3$, then $D$ contains $O_4$ as a subdigraph, a contradiction to Lemma~\ref{lemma:o4_free}. Hence $D\ind{\{u,x,y,z\}}$ contains $\ori{W_3}$ as a subdigraph. Since both $D$ and $\ori{W_3}$ are $3$-dicritical, we have $D=\ori{W_3}$.
\end{proof}

We now prove a similar result for a collection of four digraphs.
Given two digraphs $H_1$ and $H_2$, let $H_1 \Ra H_2$ denote the {\it directed join} of $H_1$ and $H_2$, that is the digraph obtained from disjoint copies of $H_1$ and $H_2$ by adding all arcs from the copy of $H_1$ to the copy of $H_2$. 
If we further add all the arcs from $H_2$ to $H_1$, we obtain the {\it bidirected join} of $H_1$ and $H_2$, denoted by $H_1 \dis H_2$. It is straightforward that $\dic(H_1 \dis H_2) = \dic(H_1) + \dic(H_2)$, see~\cite{bangjensenBSS2019}.

\begin{lemma}\label{lemma:c3_c3_free}
Let $H_1,H_2$ be two digraphs in $\{\bid{K_2},\ori{C_3}\}$ and let $D$ be a semi-complete digraph containing $H_1 \Rightarrow H_2$ as a subdigraph.  Then $D$ is 3-dicritical if and only if $D$ is exactly $\mathcal{R}(H_1,H_2)$.
\end{lemma}
\begin{proof}
    It is easy to see that $\mathcal{R}(H_1,H_2)$ is 3-dicritical. For the other direction, let us fix $H_1,H_2 \in \{\bid{K_2},\ori{C_3}\}$ and let $D$ be a 3-dicritical semi-complete digraph containing $H_1 \Rightarrow H_2$.

    Let $X= V(H_1)$ and $Y=V(H_2)$. Let us first prove that $V(D) \setminus (X \cup Y) \neq \emptyset$, so assume for a contradiction that $V(D) = X\cup Y$. 
    We claim that there exists a simple arc $uv$ from $X$ to $Y$. If this is not the case, then $D$ is exactly $H_1 \dis H_2$, so it has dichromatic number $4$, a contradiction. This simple arc $uv$ belongs to an induced directed triangle by Lemma~\ref{lemma:every_arc_in_triangle}. This directed triangle uses an arc from $Y$ to $X$, which is necessarily in a digon, a contradiction since it must be induced. Henceforth we assume that $V(D) \setminus (X \cup Y) \neq \emptyset$.

    First suppose that there exists some $v \in V(D)\setminus (X \cup Y)$ having at least one in-neighbour and one out-neighbour in both $X$ and $Y$. Since $H_1$ and $H_2$ are strongly connected, there exist four distinct vertices $x_1,x_2\in X$ and $y_1,y_2\in Y$ such that $\{x_1x_2,y_1y_2, x_1v, vx_2, y_1v,vy_2\}$ are all arcs of $D$. Then $D\ind{\{x_1,v,x_2,y_1,y_2\}}$ contains  $O_5$ as a subdigraph, a contradiction to Lemma~\ref{lemma:o5_free}. Henceforth we may assume that every vertex $v\in V(D)\setminus (X\cup Y)$ has no out-neighbour or no in-neighbour in one of $\{X,Y\}$.

    Now suppose that there exists some $v \in V(T)\setminus (X \cup Y)$ that dominates $X$. If $v$ has an out-neighbour $y$ in $Y$, then $D\ind{X \cup \{v,y\}}$ contains $O_4$ as a subdigraph if $H_1=\bid{K_2}$ and $O_5$ otherwise, a contradiction to Lemma~\ref{lemma:o4_free} or~\ref{lemma:o5_free}, respectively. Hence $v$ has no out-neighbour in $Y$. Since $D$ is semi-complete, this implies that $Y$ dominates $v$. Hence $D$ contains $\mathcal{R}(H_1,H_2)$, implying that $D$ is exactly $\mathcal{R}(H_1,H_2)$ since both $D$ and $\mathcal{R}(H_1,H_2)$ are $3$-dicritical.
    
    Henceforth we assume that for every vertex $v \in V(D)\setminus (X \cup Y)$, there exists in $D$ a simple arc from $X$ to $v$. By directional duality, there exists also a simple arc from $v$ to $Y$. Recall that every vertex $v\in V(D)\setminus (X\cup Y)$ has no out-neighbour or no in-neighbour in one of $\{X,Y\}$. We conclude on the existence of a partition $(V_1,V_2)$ of $V(D)\setminus (X \cup Y)$ such that there is no arc from $V_1$ to $X$ and there is no arc from $Y$ to $V_2$.

    By symmetry, we may assume that $V_2$ is non-empty. Let us fix $v_2\in V_2$ and $y_1\in Y$. Since $v_2y_1$ is a simple arc, by Lemma~\ref{lemma:every_arc_in_triangle}, there exists a vertex $v_1$ such that $v_2y_1v_1v_2$ is an induced directed triangle in $D$. Note that $v_1\notin Y$ since there is no arc from $Y$ to $V_2$. Also note that $v_1\notin X$ for otherwise $v_2y_1v_1v_2$ is not induced since $X$ dominates $Y$. Further note that $v_1\notin V_2$ since it is an out-neighbour of $y_1$. This implies $v_1 \in V_1$. As $v_2$ does not dominate $X$, there is some vertex in $X$, say $x_1$ that dominates $v_2$. Note that $x_1$ dominates $v_1$ by definition of $V_1$. Let $y_2$ be the unique out-neighbour of  $y_1$ in $Y$. 
    
    If $v_1$ dominates $y_2$, we obtain that $D\ind{\{x_1,v_1,v_2,y_1,y_2\}}$ contains  $O_5$ as a subdigraph, a contradiction to Lemma~\ref{lemma:o5_free}. We may hence suppose that $y_2$ dominates $v_1$. As $Y$ does not dominate $v_1$, this implies that $H_2$ is $\ori{C_3}$ and the out-neighbour $y_3$ of $y_2$ is dominated by $v_1$. Then $D\ind{\{x_1,v_1,v_2,y_2,y_3\}}$ contains  $O_5$ as a subdigraph, a contradiction to Lemma~\ref{lemma:o5_free}.
\end{proof}

    The rest of the preparatory results before the main proof of Theorem~\ref{thm:all_3_dicritical_semi_complete} aims to exclude a collection of tournaments $\mathcal{T}_8$ as induced subdigraphs and another digraph $F$ as a (not necessarily induced) subdigraph. As the proofs of these results contain several common preliminaries, we give them together. While the exact definition of $\mathcal{T}_8$ is postponed, we now give the definition of $F$.
    Let $F$ be the oriented graph with vertex set $\{u_1,\dots,u_6,x_1,x_2,x_3\}$ such that:
    \begin{itemize}
        \item $\{u_1,\dots,u_6\}$ induces a copy of $TT_6$ the unique acyclic ordering of which is exactly $u_1,\dots,u_6$, and
        \item for every $i\in[3]$, $F$ contains the arcs $u_{2i}x_i$ and $x_iu_{2i-1}$.
    \end{itemize}
    See Figure~\ref{fig:graph_F} for an illustration of $F$.
    
    \begin{figure}[ht]
        \begin{center}	
              \begin{tikzpicture}[thick,scale=1, every node/.style={transform shape}]
                \tikzset{vertex/.style = {circle,fill=black,minimum size=5pt, inner sep=0pt}}
                \tikzset{edge/.style = {->,> = latex'}}
                \foreach \i in {1,...,6}{
                    \node[vertex, label=below:$u_\i$] (u\i) at (\i,0) {};
                    \ifthenelse{\i>1}{
                        \pgfmathtruncatemacro{\j}{\i-1}
                        \foreach \k in {1,...,\j}{
                            \ifthenelse{\k<\j}{
                                \draw[edge, bend left=45] (u\k) to (u\i) {};
                            }
                            {
                                \draw[edge] (u\k) to (u\i) {};
                            }
                        }
                    }{}
          	}
               \foreach \i in {1,...,3}{
                    \pgfmathtruncatemacro{\j}{2*\i-1}
                    \pgfmathtruncatemacro{\k}{2*\i}
                    \node[vertex, label=below:$x_\i$] (x\i) at (2*\i-0.5,-0.866) {};
                    \draw[edge] (u\k) to (x\i) {};
                    \draw[edge] (x\i) to (u\j) {};
              }
              \end{tikzpicture}
          \caption{The oriented graph $F$.}
          \label{fig:graph_F}
        \end{center}
    \end{figure}

    We let $\mathcal{F}$ be the set of tournaments $T$ with vertex set $V(T) = V(F)$ and such that $A(F) \subseteq A(T)$. Note that $\mathcal{F}$ contains $2^{15}$ tournaments since $F$ has exactly $15$ pairs of non-adjacent vertices. Four of them are of special interest and we denote them by $T^1,\dots,T^4$. We give their adjacency matrices in Appendix~\ref{appendix:adjcency_matrices}.

    \begin{lemma}\label{lemma:T1_T2_T3_T4}
        None of the tournaments in $\mathcal{F} \setminus \{T^1,T^2,T^3,T^4\}$ is a subdigraph of a $3$-dicritical semi-complete digraph.
     \end{lemma}

    \begin{proof}
        For every tournament $T\in \mathcal{F}$, we check, using the code of Appendix~\ref{b1}, if it contains  $\ori{C_3} \Ra \ori{C_3}$ as a subdigraph or if it admits no $uv$-colouring for an arc $uv$. This is always the case except when $T\in \{T^1,T^2,T^3,T^4\}$. The claim then follows by Lemmas~\ref{lemma:dicol} and~\ref{lemma:c3_c3_free}.
    \end{proof}
    
    Let $F^+$ be the oriented graph obtained from $F$ by adding a vertex $u_0$ and the arcs of $\{u_0u_i \mid i\in [6]\}$. Analogously, let $F^-$ be the oriented graph obtained form $F$ by adding a vertex $u_7$ and the arcs of $\{u_iu_7 \mid i\in [6]\}$. See Figure~\ref{fig:graph_Fp_Fm} for an illustration.

    \begin{figure}[ht]
        \begin{center}	
              \begin{tikzpicture}[thick,scale=1, every node/.style={transform shape}]
                \tikzset{vertex/.style = {circle,fill=black,minimum size=5pt, inner sep=0pt}}
                \tikzset{edge/.style = {->,> = latex'}}
                \foreach \i in {0,...,6}{
                    \node[vertex, label=below:$u_\i$] (u\i) at (\i,0) {};
                    \ifthenelse{\i>0}{
                        \pgfmathtruncatemacro{\j}{\i-1}
                        \foreach \k in {0,...,\j}{
                            \ifthenelse{\k<\j}{
                                \draw[edge, bend left=45] (u\k) to (u\i) {};
                            }
                            {
                                \draw[edge] (u\k) to (u\i) {};
                            }
                        }
                    }{}
          	}
               \foreach \i in {1,...,3}{
                    \pgfmathtruncatemacro{\j}{2*\i-1}
                    \pgfmathtruncatemacro{\k}{2*\i}
                    \node[vertex, label=below:$x_\i$] (x\i) at (2*\i-0.5,-0.866) {};
                    \draw[edge] (u\k) to (x\i) {};
                    \draw[edge] (x\i) to (u\j) {};
              }

              \foreach \i in {1,...,7}{
                    \node[vertex, xshift=6.5cm, label=below:$u_\i$] (v\i) at (\i,0) {};
                    \ifthenelse{\i>1}{
                        \pgfmathtruncatemacro{\j}{\i-1}
                        \foreach \k in {1,...,\j}{
                            \ifthenelse{\k<\j}{
                                \draw[edge, bend left=45] (v\k) to (v\i) {};
                            }
                            {
                                \draw[edge] (v\k) to (v\i) {};
                            }
                        }
                    }{}
          	}
               \foreach \i in {1,...,3}{
                    \pgfmathtruncatemacro{\j}{2*\i-1}
                    \pgfmathtruncatemacro{\k}{2*\i}
                    \node[vertex, xshift=6.5cm, label=below:$x_\i$] (y\i) at (2*\i-0.5,-0.866) {};
                    \draw[edge] (v\k) to (y\i) {};
                    \draw[edge] (y\i) to (v\j) {};
              }
              \end{tikzpicture}
          \caption{The oriented graphs $F^+$ (left) and $F^-$ (right).}
          \label{fig:graph_Fp_Fm}
        \end{center}
    \end{figure}

    \begin{lemma}\label{plusminus}
        Let $D$ be a $3$-dicritical semi-complete digraph. Then $D$ does not contain a digraph in $\{F^+,F^-\}$ as a subdigraph.
      \end{lemma}
      \begin{proof}
       Observe that the digraph obtained from $F^-$ by reversing all its arcs is isomorphic to $F^+$. As the digraph obtained from a 3-dicritical, semi-complete digraph by reversing all arcs is 3-dicritical and semi-complete, it suffices to prove the statement for $F^+$. 

       In order to do so, suppose for the sake of a contradiction, that there is a $3$-dicritical semi-complete digraph $D$ containing $F^+$. By Lemma~\ref{lemma:T1_T2_T3_T4}, $D - u_0$ contains some $T'\in \{T^1,T^2,T^3,T^4\}$. Now consider the collection $\mathcal{T}$ of tournaments on $\{u_0,\ldots,u_6,v_1,v_2,v_3\}$ that have one of $T^1,T^2,T^3,T^4$ as a labelled subdigraph and in which $u_0$ dominates $\{u_1,\ldots,u_6\}$. Observe that by assumption, $D$ contains a tournament in $\mathcal{T}$ as a spanning subdigraph. Further, $\mathcal{T}$ contains exactly  $4\times 2^3 = 32$ digraphs. Using the code in Appendix~\ref{b2}, we check that each of them contains $\ori{C_3}\Ra \ori{C_3}$ or contains an arc $uv$ with no $uv$-colouring. We conclude that the same holds for $D$, a contradiction to Lemmas~\ref{lemma:dicol} or~\ref{lemma:c3_c3_free}.
    \end{proof}

   We are now ready to show that 3-dicritical semi-complete digraphs do not contain large transitive tournaments as induced subdigraphs.
    \begin{lemma}
        \label{lemma:TT8_free}
        Let $D$ be a $3$-dicritical semi-complete digraph. Then $D$ does not contain $TT_8$ as an induced subdigraph.
   \end{lemma}
    \begin{proof}
        For a contradiction, assume that $D=(V,A)$ is a $3$-dicritical semi-complete digraph containing $TT_8$ as an induced subdigraph. We will prove that $D$ contains $F^+$ or $F^-$, which is a contradiction to Lemma~\ref{plusminus}.

        Let $S\subseteq V$ be such that $D\ind{S}$ is isomorphic to $TT_8$. Let $v_1,\dots,v_8$ be the unique acyclic ordering of $S$. By Lemma~\ref{lemma:every_arc_in_triangle}, for every $i \in [7]$, there exists a vertex $x_i \in V\setminus S$ such that $v_iv_{i+1}x_iv_i$ forms an induced directed triangle $C_i$.

    Let $H$ be the graph with vertex set $V(H)=[7]$ and that contains an edge linking $i$ and $j$ if $V(C_i) \cap V(C_j) \neq \emptyset$. For any $i,j \in [7]$ with $ij \in E(H)$ and $|i-j|\geq 2$, we have $x_i=x_j$. By Lemma~\ref{lemma:intervals}, there is no set $\{i,j,k\}\subseteq [7]$ such that $x_i=x_j=x_k$. This yields that $H$ is obtained from a path on 7 vertices by adding a matching.
    We deduce from Lemma~\ref{lemma:matchpath} that there is a set $I\subseteq [7]$ with $|I| = 3$ such that the following hold:
    \begin{enumerate}[label=(\alph*)]
        \item $\{1,7\}\setminus I \neq \emptyset$,
        \item $C_i$ and $C_j$ are vertex-disjoint for all $\{i,j\}\subseteq I$.
    \end{enumerate}
    This shows that $D$ contains $F^+$ or $F^-$, yielding a contradiction to Lemma~\ref{plusminus}.
    \end{proof}

    Given an integer $k$ and a semi-complete digraph $D$, a \textit{$k$-extension} of $D$ is a semi-complete digraph on $n(D) + k$ vertices containing $D$ as an induced subdigraph. Given a set $S$ of semi-complete digraphs, a {\it $k$-extension} of $S$ is a semi-complete digraph that is a $k$-extension of some $D\in S$. 
    We are now ready to prove that no 3-dicritical semi-complete digraph contains $F$ as a subdigraph.
    
    \begin{lemma}\label{fsub}
        Let $D$ be a $3$-dicritical semi-complete digraph. Then $D$ does not contain  $F$ as a subdigraph.
    \end{lemma}
    
    \begin{proof}
     By Lemma~\ref{lemma:T1_T2_T3_T4}, it remains to show that $D$ does not contain a graph in $\{T^1,T^2,T^3,T^4\}$ as a subdigraph. Assume for a contradiction that $D$ contains at least one of $\{T^1,T^2,T^3,T^4\}$ as a subtournament.

    We use the code in Appendix~\ref{b3}. In a first part, we compute the set $\mathscr{L}$ of all semi-complete digraphs $L$ on nine vertices such that each of the following holds:
        \begin{itemize}
            \item[$(i)$]  $L$ contains some $T\in\{T^1,T^2,T^3,T^4\}$ as a subdigraph,
            
            \item[$(ii)$]  $L$ does not contain any digraph in $\{ \bid{S_4}, \bid{K_2} \Ra \bid{K_2}, \bid{K_2} \Ra \ori{C_3}, \ori{C_3} \Ra \bid{K_2}, \ori{C_3} \Ra \ori{C_3}, O_4, O_5\}$ as a subdigraph,
            
            \item[$(iii)$]  $L$ admits a $uv$-colouring for every arc $uv \in A(L)$, and
            
            \item[$(iv)$]  $L$ does not contain $TT_8$ as an induced subdigraph.
        \end{itemize}
        
        By Lemmas~\ref{lemma:dicol},~\ref{lemma:o5_free},~\ref{lemma:o4_free},~\ref{lemma:S4_free}, and~\ref{lemma:c3_c3_free}, we know that $D$ contains some $L\in \mathscr{L}$ as an induced subdigraph. In the second part of the code, we check that every $2$-extension $L'$ of $\mathscr{L}$ does not satisfy at least one of the properties $(ii)$, $(iii)$ and $(iv)$.
        
        This shows, by Lemmas~\ref{lemma:dicol},~\ref{lemma:o5_free},~\ref{lemma:o4_free} and~\ref{lemma:c3_c3_free}, that either $D\in \mathscr{L}$ or $D$ is a $1$-extension of $\mathscr{L}$. Finally, we check that every $L\in \mathscr{L}$ has dichromatic number at most two, and that every $1$-extension $L'$ satisfying $(ii)$, $(iii)$ and $(iv)$ has dichromatic number at most two. This yields a contradiction.
    \end{proof}

    We now give the definition of $\mathcal{T}_8$ and show that no digraph in $\mathcal{T}_8$ can be contained in a 3-dicritical semi-complete digraph as an induced subdigraph.
    Let $\mathcal{T}_8$ be the set of tournaments obtained from $TT_8$ by reversing exactly one arc. Observe that $TT_8$ belongs to $\mathcal{T}_8$.

    \begin{lemma}
        \label{lemma:T8_free}
        Let $D$ be a $3$-dicritical semi-complete digraph. Then $D$ does not contain any digraph in $\mathcal{T}_8$ as an induced subdigraph.
    \end{lemma}
    \begin{proof}
        Assume for a contradiction that $D$ contains some $T'\in \mathcal{T}_8$ as an induced subtournament. Let $X \subseteq V(T)$ be such that $D\ind{X}$ is isomorphic to $T'$. By definition of $\mathcal{T}_8$, let $x_1,\dots,x_8$ be an ordering of $X$ such that $D$ contains every arc $x_ix_j$ when $i<j$, except for exactly one pair $\{k,\ell\}$, $k<\ell$.

        Assume first that $k = \ell-1$. Then observe that $T'$ is isomorphic to $TT_8$, with the acyclic ordering obtained from $x_1,\dots,x_8$ by swapping $x_\ell$ and $x_{\ell-1}$. This contradicts Lemma~\ref{lemma:TT8_free}.

        Henceforth assume that $k\leq \ell-2$. If $k\geq 2$ and $\ell \leq 7$ then $D\ind{\{x_1,x_k,x_{\ell-1},x_\ell,x_8\}}$ is isomorphic to $O_5$, a contradiction to Lemma~\ref{lemma:o5_free}.
        Henceforth we assume that $k=1$ or $\ell = 8$. By directional duality, we assume without loss of generality that $k=1$.
        Let $S$ be the transitive induced subtournament of $D$ on vertices $X\setminus \{x_1,x_2,x_\ell\}$. We denote its acyclic ordering by $y_1,\dots,y_5$, which exactly corresponds to $x_3,\dots,x_{\ell-1},x_{\ell+1},\dots,x_8$. 
        By Lemma~\ref{lemma:every_arc_in_triangle}, for every $k \in [4]$, there exists a vertex $z_k$ such that $y_ky_{k+1}z_ky_k$ forms a directed triangle $C_k$. As $S$ is induced, $z_k$ must be in $V\setminus V(S)$. Moreover, $z_k \notin \{x_1,x_2,x_\ell\}$ because both $X\setminus \{x_1\}$ and $X\setminus \{x_\ell\}$ are acyclic.
        
        Let $H$ be the graph with vertex set $V(H)=[4]$ and that contains an edge linking $i$ and $j$ if $V(C_i) \cap V(C_j) \neq \emptyset$. For any $i,j\in [4]$ with $ij \in E(H)$ and $|i-j|\geq 2$, we have $z_i=z_j$. By Lemma~\ref{lemma:intervals}, there is no set $\{h,i,j\}\subseteq [4]$ such that $z_i=z_j=z_h$. This yields that $H$ is obtained from a path on 4 vertices by adding a matching containing at most 2 edges. Hence $H$ contains two non-adjacent vertices, corresponding to two disjoint directed triangles $C_i$ and $C_j$ in $D$. Together with the directed cycle $C_h = x_1x_2x_\ell$, we deduce that $D$ contains $F$ as a subdigraph. This contradicts Lemma~\ref{fsub}.
    \end{proof}
We have now proved all necessary structural properties of 3-dicritical semi-complete digraphs. The following result contains the decisive step of the proof and it requires heavy computation. 
    For every $i\in [7]$, let $\mathscr{D}_i$ be the set of semi-complete digraphs $D$ such that each of the following holds:
    \begin{itemize}
        \item the maximum acyclic set $S\subseteq V(D)$ of $D$ has size exactly $i$,
        \item for every arc $uv$ of $D$, $D$ admits a $uv$-colouring,
        \item $D$ does not contain any digraph of $\{\bid{S_4}, \bid{K_2} \Ra \bid{K_2}, \bid{K_2} \Ra \ori{C_3}, \ori{C_3} \Ra \bid{K_2}, \ori{C_3} \Ra \ori{C_3}, O_4, O_5, F\}$ as a subdigraph,
        \item $D$ does not contain any digraph of $\mathcal{T}_8$ as an induced subdigraph,
    \end{itemize}

    \begin{lemma}\label{finale}
        The  $3$-dicritical digraphs in $\bigcup_{i=1}^7 \mathscr{D}_i$ are exactly $\bid{K_3}$, $\mathcal{H}_5$, and $\mathcal{P}_7$.
    \end{lemma}
    
    \begin{proof}
        For every $i\in [7]$, we compute $\mathscr{D}_i$ by starting from the singleton $\{TT_i\}$ which is clearly the only digraph in $\mathscr{D}_i$ on at most $i$ vertices. Using the code in Appendix~\ref{b4}, we first successively compute the digraphs in $\mathscr{D}_i$ on $j\geq i$ vertices by generating every possible $1$-extension of the digraphs in $\mathscr{D}_i$ on $j-1$ vertices, and saving only the ones satisfying the conditions on $\mathscr{D}_i$. When $j$ is large enough, it turns out that the set of digraphs in $\mathscr{D}_i$ on $j$ vertices is empty, implying that $\mathscr{D}_i$ is finite.

        We then consider every digraph $D \in \mathscr{D}_i$ and check whether $D$ is $2$-dicolourable. When it is not, since it admits a $uv$-colouring for every arc $uv$, we conclude that $D$ is $3$-dicritical.
    \end{proof}

We are now ready to conclude the proof of Theorem~\ref{thm:all_3_dicritical_semi_complete}.
\begin{proof}
    By Lemmas~\ref{lemma:dicol},~\ref{lemma:o5_free},~\ref{lemma:o4_free},~\ref{lemma:S4_free},~\ref{lemma:c3_c3_free},~\ref{fsub}, and~\ref{lemma:T8_free}, we have that every 3-dicritical semi-complete digraph that is not contained in $\bigcup_{i=1}^7 \mathscr{D}_i$ is one of $\ori{W_3}$, $\mathcal{R}(\bid{K_2}, \bid{K_2})$, $\mathcal{R}(\bid{K_2},\ori{C_3})$, $\mathcal{R}(\ori{C_3}, \bid{K_2})$, and $\mathcal{R}(\ori{C_3}, \ori{C_3})$. The statement then follows directly from Lemma~\ref{finale}.
\end{proof}

\section{Maximum number of arcs in 3-dicritical digraphs}
\label{sec:bound_two_third}

This section is devoted to the proof of Theorem~\ref{thm:bound_two_third}. We need a collection of intermediate results. 
We first show that the bidirected part of a $3$-dicritical digraph is a forest unless $D$ is a bidirected odd cycle.
\begin{proposition}
    \label{prop:digon_forest}
   Let $D$ be a $3$-dicritical digraph that is not a bidirected odd cycle. Then $B(D)$ is a forest.
\end{proposition}
\begin{proof}
 Assume for a contradiction that $B(D)$ is not a forest. Then it contains a cycle $C=u_1u_2\dots u_{p} u_1$. Let $\bid{C}$ be the bidirected cycle in $D$ corresponding to $C$. The cycle $C$ cannot be odd, for otherwise $\bid{C}$ would be a bidirected odd cycle, and $D=\bid{C}$ because a bidirected odd cycle is $3$-dicritical, a contradiction.
 Hence $C$ is an even cycle. 
 By Lemma~\ref{lemma:dicol}, there exists a $2$-dicolouring $\phi$ of $D\setminus \{u_1u_p\}$. Necessarily, $u_1$ and $u_p$ are coloured differently because there is a bidirected path of odd length between $u_1$ and $u_p$. Thus $\phi$ is a $2$-dicolouring of $D$, a contradiction.
\end{proof}

    For the remainder of this section we need a few specific definitions. Let $T$ be a tree and $V_3(T)$ be the set of vertices of degree at least $3$ in $T$. Two vertices $u,v \in V(T)$ form an {\it odd pair} if they are non-adjacent and $\dist_T(u,v)$ is odd, where $\dist_T(u,v)$ denotes the length of the unique path between $u$ and $v$ in $T$. The set of odd pairs of $T$ is denoted by $\OP(T)$ and its cardinality is denoted by $\op(T)$. We finally define the {\it dearth} of $T$ as follows:
    \[
        \dearth(T) = \sum_{v\in V_3(T)}\frac{1}{6}d(v)(d(v)-1) + \op(T).
    \]
    
    We first prove that the dearth of a tree is always at least a fraction of its order.

    \begin{lemma}
        \label{lemma:dearth_1_3_n}
        Let $T$ be a tree on $n$ vertices for some positive integer $n$. Then $\dearth(T) \geq \frac{1}{3}n - 1$.
    \end{lemma}
    \begin{proof}
        For the sake of a contradiction, suppose that $T$  is a counterexample to the statement whose number of vertices is minimum. Clearly, we have $n\geq 4$. The following claim excludes a collection of simple structures of $T$.
        \begin{claim}\label{clm:cheminetoile}
            $T$ is neither a path nor a tree.
        \end{claim}
        \begin{proofclaim}
        The statement follows from the following simple case distinction.
        \begin{description}
            \item[\textbf{Case 1:}]\textit{$T$ is a path of even length.}
            
            For every odd $i\in \{1,\dots,n-3\}$, as $n\geq 4$, there are exactly $i$ distinct pairs of vertices at distance exactly $n-i$ in $T$. Hence $\dearth(T) \geq \op(T) = \sum_{i=1}^{\frac{n-2}{2}}(2i-1) = \left( \frac{n-2}{2} \right)^2\geq \frac{1}{3}n-1$.
            
            \item[\textbf{Case 2:}]\textit{$T$ is a path of odd length.}
            
            For every odd $i\in \{1,\dots,n-2\}$, as $n \geq 4$, there are exactly $i$ distinct pairs of vertices at distance exactly $n-i$ in $T$. Hence $\dearth(T) \geq \op(T) = \sum_{i=1}^{\frac{n-3}{2}}2i = \left( \frac{n-3}{2} \right) \left( \frac{n-1}{2} \right) \geq \frac{1}{3}n-1$.
            
            \item[\textbf{Case 3:}]\textit{$T$ is a star on $n\geq 4$ vertices.}
            
            As $n\geq 4$, we obtain that $\dearth(T)$ is exactly $\frac{1}{6}(n-1)(n-2)$, and so  $\dearth(T)\geq \frac{1}{3}n-1$.
             \end{description}
             In either case, we obtain a contradiction to the choice of $T$.
        \end{proofclaim}

            By Claim \ref{clm:cheminetoile}, we obtain that $T$ is neither a path nor a star. In particular, it follows that $T$ contains an edge $uv$ such that $d_T(u)\geq 2$ and $d_T(v) \geq 3$.
            Let $v_1,\dots,v_r$ be the neighbours of $v$ in $T$, where $v_1 = u$ and $r=d_T(v)\geq 3$. 
            For each $i\in[r]$, let $T_i$ be the component of $T-v$ containing $v_i$. By the choice of $T$, we have $\dearth(T_i) \geq \frac{1}{3}n(T_i)-1$.
            Since the $T_i$s are pairwise disjoint and none of them contains $v$, and because $u$ has a neighbour in $T_1$ at distance exactly $3$ from $v_2,\dots,v_r$ we obtain:
            \begin{align*}
                \dearth(T) &\geq \sum_{i=1}^{r} \dearth(T_i) + \frac{1}{6}r(r-1) + (r-1)\\
                &\geq \frac{1}{3}(n(T)-1) - r + \frac{1}{6}r(r-1) + (r-1)\\
                &\geq \frac{1}{3}n(T) - 1,
            \end{align*}
            where in the last inequality we used $r\geq 3$. This contradicts the choice of $T$.
        \end{proof}

    \begin{lemma}
        \label{lemma:dearth}
        Let $D$ be a $3$-dicritical digraph distinct from $\{\bid{K_3},\ori{W_3}\}$ and $\bid{T}$ be a bidirected tree contained in $D$. Then we have  
        \[
            |\{ \{u,v\} \subseteq V(T) \mid \{uv,vu\}\cap A(D) = \emptyset \}| \geq \dearth(T).
        \]
    \end{lemma}
    \begin{proof}
        Set $\mathcal{O} = \{ \{u,v\} \subseteq V(T) \mid \{uv,vu\}\cap A(D) = \emptyset \}$. For every vertex $v\in V_3(T)$, let $\mathcal{O}_v = \mathcal{O} \cap (N_T(v) \times N_T(v))$. Finally let $\mathcal{O}_{\mathrm{odd}} = \mathcal{O} \cap \OP(T)$. 

        Let us first show that these sets are pairwise disjoint. Let $u,v\in V_3(T)$ be two vertices of degree at least $3$ in $T$. Since $T$ is a tree, we have that $N_T(u)\cap N_T(v)$ contains at most one vertex, implying that $\mathcal{O}_u \cap \mathcal{O}_v = \emptyset$. Also note that vertices in $N_T(v)$ are at distance exactly $2$ from each other, so $\mathcal{O}_v \cap \mathcal{O}_{\mathrm{odd}} = \emptyset$. This implies
        \[
            |\mathcal{O}| \geq \sum_{v\in V_3(T)}|\mathcal{O}_v| + |\mathcal{O}_{\mathrm{odd}}|.
        \]
        Hence it is sufficient to prove $|\mathcal{O}_v| \geq \frac{1}{6}d_T(v)(d_T(v)-1)$ for every $v\in V_3(T)$ and $\mathcal{O}_{\mathrm{odd}} = \OP(T)$ to prove Lemma \ref{lemma:dearth}.
        
        Let $v\in V_3(T)$ and $u,x,z$ be three distinct vertices in $N_T(v)$. We claim that $D\ind{\{u,x,z\}}$ contains at most two arcs. If this is not the case, then $D\ind{\{u,x,z\}}$ contains a digon, a directed triangle or a transitive tournament on three vertices. This implies that $D\ind{\{u,x,z,v\}}$ contains $\bid{K_3}$, $\ori{W_3}$ or $O_4$. By Theorem \ref{thm:all_3_dicritical_semi_complete} and Lemma \ref{lemma:o4_free}, in each case, we obtain a contradiction to the choice of $D$.
        Since this holds for every choice of three distinct vertices in $N_T(v)$ and each pair of vertices in $N_T(v)$ is contained in $d_T(v)-2$ triples, we deduce the following inequality 
        \[
            m(D\ind{N_T(v)}) \cdot (d_T(v)-2) = \sum_{\substack{X\subseteq N_T(v),\\|X|=3}} m(D\ind{X}) \leq 2 \cdot \binom{d_T(v)}{3},
        \]
        implying that $m(D\ind{N_T(v)}) \leq \frac{1}{3}d_T(v)(d_T(v)-1)$.
        Therefore, we obtain $|\mathcal{O}_v| = \binom{d_T(v)}{2} - m(D\ind{N_T(v)}) \geq \frac{1}{6}d_T(v)(d_T(v)-1)$ as desired.
        
        To show $\mathcal{O}_{\mathrm{odd}} = \OP(T)$, it is sufficient to show that if $\{u,v\}$ is an odd pair then $\{uv,vu\}\cap A(D) = \emptyset$. Assume this is not the case, then by Lemma~\ref{lemma:dicol} $D' = D \setminus \{uv,vu\}$ admits a $2$-dicolouring $\phi$ in which $\phi(u) = \phi(v)$, a contradiction since $u$ and $v$ are connected by a bidirected odd path in $D'$. This shows the claim.
    \end{proof}

We are now ready to prove Theorem~\ref{thm:bound_two_third} that we first restate here for convenience.

\thmboundtwothird*

\begin{proof}
    Let $D$ be such a digraph. If $D$ is a bidirected odd cycle, we have $n \geq 5$ and hence $D$ has $2n \leq \binom{n}{2} + \frac{2}{3}n$ arcs, so the result trivially holds. Henceforth assume $D$ is not a bidirected cycle. Then, by Proposition~\ref{prop:digon_forest}, $B(D)$ is a forest. Let $T_1,\dots,T_s$ be the connected components of $B(D)$. For every $i\in [s]$, the number of digons in $D\ind{V(T_i)}$ is exactly $n(T_i)-1$, whereas the number of pairs of non-adjacent vertices is at least $\dearth(T_i)$ by Claim~\ref{lemma:dearth}. Hence, since there is no digon between the $T_i$s, we obtain
    \begin{align*}
        m(D) &\leq \binom{n}{2} + \sum_{i=1}^s  \big((n(T_i)-1) - \dearth(T_i) \big)\\
        &\leq \binom{n}{2} + \sum_{i=1}^s  \left((n(T_i)-1) - (\frac{1}{3}n(T_i) -1) \right) \text{~~~by Claim~\ref{lemma:dearth_1_3_n}}\\
        &= \binom{n}{2} + \frac{2}{3}n,
    \end{align*} 
    which concludes the proof.
\end{proof}

\section{Conclusion}\label{conc}

In this paper, we showed that the number of 3-dicritical semi-complete digraphs is finite and with a computer-assisted proof, we gave a full characterization of them. This result seems to be only the tip of an iceberg, and natural generalizations in several directions can be considered.

\medskip

First, the conjecture of Hoshino and Kawarabayashi on the maximum density of $3$-dicritical oriented graphs remains widely open. 

\medskip 

We believe that almost all 3-dicritical digraphs are sparser than tournaments.
We thus propose the following conjecture which would imply Theorem~\ref{thm:finite_3_dicritical_semi_complete} and asymptotically improve on Theorem~\ref{thm:bound_two_third}.

\begin{conjecture}\label{conj:dense}
    There is only a finite number of 3-dicritical digraphs $D$ on $n$ vertices that satisfy $m(D)\geq \binom{n}{2}$.
\end{conjecture}

\medskip

It is an interesting challenge to generalize the results obtained in this article to $k \geq 4$. In particular, we would be interested in a confirmation of the following statement.

\begin{conjecture}
    \label{conj:finiteness_k}
    For every $k \geq 4$, there is only a finite number of $k$-dicritical semi-complete digraphs.
\end{conjecture}

\medskip

Finally, it is also natural to consider a different notion of criticality. A digraph $D$ is called {\it 3-vertex-dicritical} if $D$ is not 3-dicolourable, but $D-v$ is for all $v \in V(D)$. Observe that every 3-dicritical digraph is 3-vertex-dicritical, but the converse is not necessarily true. One can hence wonder whether an analogue of Theorem~\ref{thm:finite_3_dicritical_semi_complete} is true for 3-vertex-dicritical digraphs. However, this turns out not to be the case. 
Chen et al. proved in~\cite{chenSIAM3} that it is NP-hard to decide whether a given tournament is 2-dicolourable. An infinite collection of 3-vertex-dicritical tournaments can easily be derived from their proof.

\bibliographystyle{abbrv}
\bibliography{refs}

\newpage

\appendix

\appendix

\section{The tournaments \texorpdfstring{$T^1,\ldots,T^4$}{T1,...,T4}}
\label{appendix:adjcency_matrices}

We give the adjacency matrices of $T^1, T^2, T^3$ and $T^4$.

\[
   T^1 :~~~ 
     \begin{blockarray}{cccccccccc}
         & u_1 & u_2 & u_3 & u_4 & u_5 & u_6 & x_1 & x_2 & x_3 \\
       \begin{block}{c[ccccccccc]}
            u_1 & 0 & 1 & 1 & 1 & 1 & 1 & 0 & 1 & 1\\
            u_2 & 0 & 0 & 1 & 1 & 1 & 1 & 1 & 1 & 1\\
            u_3 & 0 & 0 & 0 & 1 & 1 & 1 & 0 & 0 & 0\\
            u_4 & 0 & 0 & 0 & 0 & 1 & 1 & 0 & 1 & 0\\
            u_5 & 0 & 0 & 0 & 0 & 0 & 1 & 0 & 1 & 0\\
            u_6 & 0 & 0 & 0 & 0 & 0 & 0 & 1 & 1 & 1\\
            x_1 & 1 & 0 & 1 & 1 & 1 & 0 & 0 & 0 & 1\\
            x_2 & 0 & 0 & 1 & 0 & 0 & 0 & 1 & 0 & 1\\
            x_3 & 0 & 0 & 1 & 1 & 1 & 0 & 0 & 0 & 0\\
       \end{block}
     \end{blockarray}
 \]
 \[
   T^2 :~~~ 
     \begin{blockarray}{cccccccccc}
         & u_1 & u_2 & u_3 & u_4 & u_5 & u_6 & x_1 & x_2 & x_3 \\
       \begin{block}{c[ccccccccc]}
            u_1 & 0 & 1 & 1 & 1 & 1 & 1 & 0 & 1 & 1\\
            u_2 & 0 & 0 & 1 & 1 & 1 & 1 & 1 & 1 & 1\\
            u_3 & 0 & 0 & 0 & 1 & 1 & 1 & 1 & 0 & 0\\
            u_4 & 0 & 0 & 0 & 0 & 1 & 1 & 0 & 1 & 0\\
            u_5 & 0 & 0 & 0 & 0 & 0 & 1 & 0 & 1 & 0\\
            u_6 & 0 & 0 & 0 & 0 & 0 & 0 & 0 & 1 & 1\\
            x_1 & 1 & 0 & 0 & 1 & 1 & 1 & 0 & 1 & 0\\
            x_2 & 0 & 0 & 1 & 0 & 0 & 0 & 0 & 0 & 1\\
            x_3 & 0 & 0 & 1 & 1 & 1 & 0 & 1 & 0 & 0\\
       \end{block}
     \end{blockarray}
 \]
 
 \[
   T^3 :~~~ 
     \begin{blockarray}{cccccccccc}
         & u_1 & u_2 & u_3 & u_4 & u_5 & u_6 & x_1 & x_2 & x_3 \\
       \begin{block}{c[ccccccccc]}
            u_1 & 0 & 1 & 1 & 1 & 1 & 1 & 0 & 0 & 0\\
            u_2 & 0 & 0 & 1 & 1 & 1 & 1 & 1 & 0 & 1\\
            u_3 & 0 & 0 & 0 & 1 & 1 & 1 & 1 & 0 & 1\\
            u_4 & 0 & 0 & 0 & 0 & 1 & 1 & 1 & 1 & 1\\
            u_5 & 0 & 0 & 0 & 0 & 0 & 1 & 0 & 0 & 0\\
            u_6 & 0 & 0 & 0 & 0 & 0 & 0 & 0 & 0 & 1\\
            x_1 & 1 & 0 & 0 & 0 & 1 & 1 & 0 & 1 & 1\\
            x_2 & 1 & 1 & 1 & 0 & 1 & 1 & 0 & 0 & 0\\
            x_3 & 1 & 0 & 0 & 0 & 1 & 0 & 0 & 1 & 0\\
       \end{block}
     \end{blockarray}
 \]
 
 \[
   T^4 :~~~ 
     \begin{blockarray}{cccccccccc}
         & u_1 & u_2 & u_3 & u_4 & u_5 & u_6 & x_1 & x_2 & x_3 \\
       \begin{block}{c[ccccccccc]}
            u_1 & 0 & 1 & 1 & 1 & 1 & 1 & 0 & 0 & 1\\
            u_2 & 0 & 0 & 1 & 1 & 1 & 1 & 1 & 0 & 1\\
            u_3 & 0 & 0 & 0 & 1 & 1 & 1 & 1 & 0 & 1\\
            u_4 & 0 & 0 & 0 & 0 & 1 & 1 & 1 & 1 & 0\\
            u_5 & 0 & 0 & 0 & 0 & 0 & 1 & 0 & 0 & 0\\
            u_6 & 0 & 0 & 0 & 0 & 0 & 0 & 0 & 0 & 1\\
            x_1 & 1 & 0 & 0 & 0 & 1 & 1 & 0 & 1 & 0\\
            x_2 & 1 & 1 & 1 & 0 & 1 & 1 & 0 & 0 & 1\\
            x_3 & 0 & 0 & 0 & 1 & 1 & 0 & 1 & 0 & 0\\
       \end{block}
     \end{blockarray}
 \]

\section{Code used in the proof of Theorem~\ref{thm:all_3_dicritical_semi_complete}}
This appendix contains the code used in the proof of Theorem~\ref{thm:all_3_dicritical_semi_complete}. In Appendix~\ref{b0}, we give a collection of useful subroutines we use in the main part of the code. In Appendices~\ref{b1},~\ref{b2},~\ref{b3}, and~\ref{b4}, we give the code use in the proofs of Lemmas~\ref{lemma:T1_T2_T3_T4} and~\ref{plusminus}, and Lemmas~\ref{fsub} and~\ref{finale}, respectively.
\setcounter{subsection}{-1}
\subsection{Preliminaries for the code}\label{b0}

In the following code, we give a collection of subroutines we use in our code.

\begin{lstlisting}[language=Python]
# The following function displays a progress bar
def printProgressBar (iteration, total):
    percent = ("{0:.1f}").format(100 * (iteration / float(total)))
    filledLength = int(50 * iteration // total)
    bar = "#" * filledLength + "-" * (50 - filledLength)
    print(f"\rProcess: |{bar}| {percent}% Complete", end = "\r")
    # Print New Line on Complete
    if iteration == total: 
        print()

# k,n : integers such that k < 3**n
# Returns: the decomposition of k in base 3 of length n
def ternary(k,n):
    b = 3**(n-1)
    res=""
    for i in range(n):
        if(k >= 2*b):
            k -= 2*b
            res = res + "2"
        elif(k >= b):
            k -= b
            res = res + "1"
        else:
            res = res + "0"
        b /= 3
    return res

# d: DiGraph
# u: vertex
# v: vertex
# Returns: True if and only if d contains a directed path from u to v
def contains_directed_path(d,u,v):
    to_be_treated = [u]
    i=0
    while(len(to_be_treated) != i):
        x = to_be_treated[i]
        if (x==v):
            return True
        for y in d.neighbors_out(x):
            if (not y in to_be_treated):
                to_be_treated.append(y)
        i+=1
    return False

# d: DiGraph
# u: vertex of d
# v: vertex of d
# current_colouring: partial 2-dicolouring with colours {0,1} of d such that current_colouring[u] = current_colouring[v] = 0
# Returns: True if and only if current_colouring can be extended into a 2-dicolouring of d with no monochromatic directed path from u to v.
def can_be_subgraph_of_3_dicritical_aux(d,u,v, current_colouring):
    #build the colour classes
    colours = {}
    colours[0] = []
    colours[1] = []
    for (x,i) in current_colouring.items():
        colours[i].append(x)
    #check whether both colour classes are acyclic
    for i in range(2):
        d_i = d.subgraph(colours[i])
        if(not d_i.is_directed_acyclic()):
            return False
    #check whether there is a monochromatic directed path from u to v
    d_0 = d.subgraph(colours[0])
    if(contains_directed_path(d_0,u,v)):
        return False
    #check whether current_colouring is partial
    if(len(current_colouring) == d.order()):
        return True
    else:
        #find a vertex x that is not coloured yet
        x = len(current_colouring)
        while(x in current_colouring):
            x-=1
        #check recursively whether current_colouring can be extended to x
        for i in range(2):
            current_colouring[x] = i
            if(can_be_subgraph_of_3_dicritical_aux(d,u,v,current_colouring)):
                return True
            current_colouring.pop(x, None)
        return False

# d: DiGraph
# forbidden_subtournaments: list of DiGraphs
# Returns: True if and only if d is {forbidden_subdigraphs}-free, {forbidden_induced_subdigraphs}-free  and, for every arc (u,v) of d, d admits a uv-colouring.
def can_be_subgraph_of_3_dicritical(d, forbidden_subdigraphs, forbidden_induced_subdigraphs):
    #check whether d contains a forbidden subgraph
    for T in forbidden_subdigraphs:
        if(d.subgraph_search(T, False) != None):
            return False
    #check whether d contains a forbidden induced subgraph
    for T in forbidden_induced_subdigraphs:
        if(d.subgraph_search(T, True) != None):
            return False
    #check for every arc uv if d admits a uv-colouring.
    for e in d.edges():
        d_aux = DiGraph(len(d.vertices()))
        d_aux.add_edges(d.edges())
        d_aux.delete_edge(e)
        current_colouring = {}
        current_colouring[e[0]] = 0
        current_colouring[e[1]] = 0
        if(not can_be_subgraph_of_3_dicritical_aux(d_aux,e[0],e[1], current_colouring)):
            return False
    return True

# d: DiGraph
# Returns: True if and only if d is 2-dicolourable
def is_two_dicolourable(d):
    n = d.order()
    for bipartition in range(2**n):
        #build the binary word corresponding to the bipartition
        binary = bin(bipartition)[2:]
        while(len(binary)<(n)):
            binary = "0" + binary
        #build the bipartition
        V1 = []
        V2 = []
        for v in range(n):
            if(binary[v] == '0'):
                V1.append(v)
            else:
                V2.append(v)
        #check whether (V1,V2) is actually a dicolouring
        d1 = d.subgraph(V1)
        d2 = d.subgraph(V2)
        if(d1.is_directed_acyclic() and d2.is_directed_acyclic()):
            return True
    return False

#C3_C3 is the digraph made of two disjoint directed triangles, the vertices of one dominating the vertices of the other
C3_C3 = DiGraph(6)
for i in range(3):
    C3_C3.add_edge(i,(i+1)%3)
    C3_C3.add_edge(i+3,((i+1)%3)+3)
    for j in range(3,6):
        C3_C3.add_edge(i,j)

#F is the digraph on nine vertices made of a TT6 u1,...,u6 and the arcs of the directed triangles u1u2x1u1, u3u4x2u3, u5u6x3u5. 
F = DiGraph(9)
for i in range(6):
    for j in range(i):
        F.add_edge(j,i)
F.add_edge(6,0)
F.add_edge(1,6)
F.add_edge(7,2)
F.add_edge(3,7)
F.add_edge(8,4)
F.add_edge(5,8)

#TT8 is the transitive tournament on 8 vertices
TT8 = DiGraph(8)
for i in range(8):
    for j in range(i):
        TT8.add_edge(j,i)

#reversed_TT8 is the set of tournaments, up to isomorphism, obtained from TT8 by reversing exactly one arc
reversed_TT8 = []
for e in TT8.edges():
    rev = DiGraph(8)
    rev.add_edges(TT8.edges())
    rev.delete_edge(e)
    rev.add_edge(e[1],e[0])
    check = True
    for T in reversed_TT8:
        check = check and (not T.is_isomorphic(rev))
    if(check):
        reversed_TT8.append(rev)

#K2 is the complete digraph on 2 vertices
K2 = DiGraph(2)
K2.add_edge(0,1)
K2.add_edge(1,0)

#S4 is the bidirected star on 4 vertices
S4 = DiGraph(4)
for i in range(1,4):
    S4.add_edge(i,0)
    S4.add_edge(0,i)

#C3_K2 is the digraph with a directed triangle dominating a digon. 
C3_K2 = DiGraph(5)
for i in range(3):
    C3_K2.add_edge(i,(i+1)%3)
    for j in range(3,5):
        C3_K2.add_edge(i,j)
C3_K2.add_edge(3,4)
C3_K2.add_edge(4,3)

#K2_C3 is the digraph with a digon dominating a directed triangle
K2_C3 = DiGraph(5)
for i in range(3):
    K2_C3.add_edge(i,(i+1)%3)
    for j in range(3,5):
        K2_C3.add_edge(j,i)
K2_C3.add_edge(3,4)
K2_C3.add_edge(4,3)

#K2_K2 is the digraph with a digon dominating a digon
K2_K2 = DiGraph(4)
for i in range(2):
    K2_K2.add_edge(i,(i+1)%2)
    K2_K2.add_edge(i+2,((i+1)%2)+2)
    for j in range(2,4):
        K2_K2.add_edge(i,j)

#O4 and O5 are the obstructions described in the paper.
O4 = DiGraph(4)
for i in range(1,4):
    O4.add_edge(0,i)
for i in range(1,3):
    O4.add_edge(i,3)
O4.add_edge(1,2)
O4.add_edge(2,1)

O5 = DiGraph(5)
for i in range(1,5):
    O5.add_edge(0,i)
for i in range(1,4):
    O5.add_edge(i,4)
O5.add_edge(1,2)
O5.add_edge(2,3)
O5.add_edge(3,1)
\end{lstlisting}

\subsection{The proof of Lemma~\ref{lemma:T1_T2_T3_T4}}\label{b1}

We here give the code used in the proof of Lemma~\ref{lemma:T1_T2_T3_T4}.

\begin{lstlisting}[language=Python]
load("tools.sage")

# binary_code: a string of fifteen characters '0' and '1'
# Returns: a tournament of \mathcal{F}. The orientations of the fifteen non-forced arcs correspond to the characters of binary_code.
def digraph_blowup_TT3(binary_code):
    iterator_binary_code = iter(binary_code)
    d = DiGraph(9)
    #the vertices 0,...,8 correspond respectively to u_1,...,u_6,x_1,x_2,x_3
    
    #add the arcs of the TT_6
    for i in range(6):
        for j in range(i):
            d.add_edge(j,i)
            
    #add the arcs of the directed triangles
    for i in range(3):
        d.add_edge(6+i, 2*i)
        d.add_edge(2*i+1, 6+i)
    
    missing_edges=[(6,2),(6,3),(6,4),(6,5),(6,7),(6,8),(7,0),(7,1),(7,4),(7,5),(7,8),(8,0),(8,1),(8,2),(8,3)]
    #we orient the missing_edges according to binary_code
    for e in missing_edges:
        if(next(iterator_binary_code) == '0'):
            d.add_edge(e[0],e[1])
        else:
            d.add_edge(e[1],e[0])
    return d

print("Computing all possible candidates of \mathcal{F} for being a subtournament of a 3-dicritical semi-complete digraph...")
list_candidates = []
list_forbidden_induced_subdigraphs = [C3_C3]

#print progress bar
printProgressBar(0, 2**15)
for i in range(2**15):
    binary_value = bin(i)[2:]
    while(len(binary_value)<15):
        binary_value = '0' + binary_value
    d = digraph_blowup_TT3(binary_value)
    if(can_be_subgraph_of_3_dicritical(d,[],list_forbidden_induced_subdigraphs)):
        list_candidates.append(d)
    #update progress bar
    printProgressBar(i + 1, 2**15)

print("Number of candidates: ",len(list_candidates),".")
for i in range(len(list_candidates)):
    print("Candidate ",i+1,": ")
    list_candidates[i].export_to_file("T"+str(i+1)+".pajek")
    print(list_candidates[i].adjacency_matrix())
\end{lstlisting}

Running this code produces the following output after roughly 2 minutes of execution on a standard desktop computer:

\begin{lstlisting}
Computing all possible candidates of \mathcal{F} for being a subtournament of a 3-dicritical semi-complete digraph...
Process: |##################################################| 100.0% Complete
Number of candidates:  4 .
Candidate  1 : 
[0 1 1 1 1 1 0 1 1]
[0 0 1 1 1 1 1 1 1]
[0 0 0 1 1 1 0 0 0]
[0 0 0 0 1 1 0 1 0]
[0 0 0 0 0 1 0 1 0]
[0 0 0 0 0 0 1 1 1]
[1 0 1 1 1 0 0 0 1]
[0 0 1 0 0 0 1 0 1]
[0 0 1 1 1 0 0 0 0]
Candidate  2 : 
[0 1 1 1 1 1 0 1 1]
[0 0 1 1 1 1 1 1 1]
[0 0 0 1 1 1 1 0 0]
[0 0 0 0 1 1 0 1 0]
[0 0 0 0 0 1 0 1 0]
[0 0 0 0 0 0 0 1 1]
[1 0 0 1 1 1 0 1 0]
[0 0 1 0 0 0 0 0 1]
[0 0 1 1 1 0 1 0 0]
Candidate  3 : 
[0 1 1 1 1 1 0 0 0]
[0 0 1 1 1 1 1 0 1]
[0 0 0 1 1 1 1 0 1]
[0 0 0 0 1 1 1 1 1]
[0 0 0 0 0 1 0 0 0]
[0 0 0 0 0 0 0 0 1]
[1 0 0 0 1 1 0 1 1]
[1 1 1 0 1 1 0 0 0]
[1 0 0 0 1 0 0 1 0]
Candidate  4 : 
[0 1 1 1 1 1 0 0 1]
[0 0 1 1 1 1 1 0 1]
[0 0 0 1 1 1 1 0 1]
[0 0 0 0 1 1 1 1 0]
[0 0 0 0 0 1 0 0 0]
[0 0 0 0 0 0 0 0 1]
[1 0 0 0 1 1 0 1 0]
[1 1 1 0 1 1 0 0 1]
[0 0 0 1 1 0 1 0 0]
\end{lstlisting}

The graphs in the output are exactly $T_1,T_2,T_3$, and $T_4$.

\subsection{The proof of Lemma~\ref{plusminus}}\label{b2}

We here give the code used in the proof of Lemma~\ref{plusminus}.

\begin{lstlisting}[language=python]
import networkx
load("tools.sage")

list_candidates = []
list_forbidden_induced_subdigraphs = [C3_C3]

#import T1, T2, T3 and T4
for i in range(1,5):
    candidate = DiGraph(9)
    nx = networkx.read_pajek("T"+str(i)+".pajek")
    for e in nx.edges():
        candidate.add_edge(int(e[0]),int(e[1]))
    list_candidates.append(candidate)

print("We start from the ", len(list_candidates), " candidates on 9 vertices.")

#We want to prove that a 3-dicritical semi-complete digraph does not contain a digraph in {F+,F-}. By directional duality, it is sufficient to prove that it does not contain F+.
#For each candidate computed above, we try to add a new vertex that dominates the transitive tournament, and then we build every possible orientation between this vertex and the three other vertices. 
print("Computing for F+...")
list_Fp = []
printProgressBar(0, 8)
for orientation in range(2**3):
    binary = bin(orientation)[2:]
    while(len(binary)<3):
        binary = '0' + binary
    for T9 in list_candidates:
        iterator = iter(binary)
        T10 = DiGraph(10)
        T10.add_edges(T9.edges())
        for v in range(6):
            T10.add_edge(9,v)
        for v in range(6,9):
            if(next(iterator) == '0'):
                T10.add_edge(v,9)
            else:
                T10.add_edge(9,v)
        check = can_be_subgraph_of_3_dicritical(T10,[],list_forbidden_induced_subdigraphs)
        if(check):
            list_Fm.append(T2)
    printProgressBar(orientation+1, 8)

print("Number of 1-extensions of {T1,T2,T3,T4} containing F+: ",len(list_Fp))
\end{lstlisting}

Running this code produces the following output after roughly 1 second of execution on a standard desktop computer:

\begin{lstlisting}
We start from the  4  candidates on 9 vertices.
Computing for F+...
Process: |##################################################| 100.0% Complete
Number of 1-extensions of {T1,T2,T3,T4} containing F+:  0
\end{lstlisting}

\subsection{The proof of Lemma~\ref{fsub}}\label{b3}

We here give the code used in the proof of Lemma~\ref{fsub}.

\begin{lstlisting}[language=python]
import networkx
load("tools.sage")

all_candidates = []
current_candidates = []
next_candidates = []

list_forbidden_subdigraphs = [S4, K2_K2, O4, O5, K2_C3, C3_K2, C3_C3]
list_forbidden_induced_subdigraphs = [TT8]

#import the candidates T1, ..., T4 on 9 vertices:
for i in range(1,5):
    candidate = DiGraph(9)
    nx = networkx.read_pajek("T"+str(i)+".pajek")
    for e in nx.edges():
        candidate.add_edge(int(e[0]),int(e[1]))
    current_candidates.append(candidate)
print("We start from the tournaments {T^1,T^2,T^3,T^4} on 9 vertices, and look for every possible completion of them that is potentially a subdigraph of a larger 3-dicritical semi-complete digraphs.\n")


all_candidates.extend(current_candidates)
#completions of T1, ..., T4
while(len(current_candidates)>0):
    for old_D in current_candidates:
        for e in old_D.edges():
            #we try to complete old_D by replacing e by a digon. It actually makes sense only if e is not already in a digon.
            if(not (e[1],e[0],None) in old_D.edges()):
                new_D = DiGraph(9)
                new_D.add_edges(old_D.edges())
                new_D.add_edge(e[1],e[0])
                #we check whether this completion of old_D is potentially a subdigraph of a larger 3-dicritical semi-complete digraph.
                check = can_be_subgraph_of_3_dicritical(new_D,list_forbidden_subdigraphs, list_forbidden_induced_subdigraphs)
                for D in next_candidates:
                    check = check and (not D.is_isomorphic(new_D))
                if(check):
                    next_candidates.append(new_D)
    all_candidates.extend(next_candidates)
    current_candidates = next_candidates
    next_candidates=[]

print("-----------------------------------------")
print("There are",len(all_candidates),"possible completions (up to isomorphism) of {T^1,T^2,T^3,T^4} that are potentially subdigraphs of a larger 3-dicritical semi-complete digraphs.\n")
                

count_dic_3 = 0
for D in all_candidates:
    if(not is_two_dicolourable(D)):
        count_dic_3 += 1
print(count_dic_3, " of them have dichromatic number at least 3. In particular,",count_dic_3,"of them are 3-dicritical.\n")

current_candidates = all_candidates
next_candidates = []
for n in range(10,12):
    #computes the extensions on n vertices of {T1,T2,T3,T4}   
    print("-----------------------------------------")
    print("Computing "+str(n-9)+"-extensions of the candidates on 9 vertices that are potentially subtournaments of 3-dicritical tournaments (up to isomorphism).")
    printProgressBar(0, 3**(n-1))

    for orientation in range(3**(n-1)):
        ternary_code = ternary(orientation,n-1)
        #build every 1-extension of current_candidates
        for old_D in current_candidates:
            new_D = DiGraph(n)
            new_D.add_edges(old_D.edges())
            for v in range(n-1):
                if(ternary_code[v] == '0'):
                    new_D.add_edge(v,n-1)
                elif(ternary_code[v]=='1'):
                    new_D.add_edge(n-1,v)
                else:
                    new_D.add_edge(v,n-1)
                    new_D.add_edge(n-1,v)
            check = can_be_subgraph_of_3_dicritical(new_D,list_forbidden_subdigraphs,list_forbidden_induced_subdigraphs)
            for D in next_candidates:
                check = check and (not new_D.is_isomorphic(D))
            if(check):
                next_candidates.append(new_D)
        printProgressBar(orientation+1, 3**(n-1))

    print("Number of ", n-9,"-extensions up to isomorphism: ",len(next_candidates))
    #check if one of the candidates has dichromatic number at least 3.
    count_dic_3 = 0
    for D in next_candidates:
        if(not is_two_dicolourable(D)):
            count_dic_3 += 1
    print(count_dic_3, " of them have dichromatic number at least 3. In particular,", count_dic_3,"of them are 3-dicritical.\n")
    current_candidates = next_candidates
    next_candidates = []
\end{lstlisting}

Running this code produces the following output after roughly 12 minutes of execution on a standard desktop computer:

\begin{lstlisting}
We start from the tournaments {T^1,T^2,T^3,T^4} on 9 vertices, and look for every possible completion of them that is potentially a subdigraph of a larger 3-dicritical semi-complete digraphs.

-----------------------------------------
There are 14 possible completions (up to isomorphism) of {T^1,T^2,T^3,T^4} that are potentially subdigraphs of a larger 3-dicritical semi-complete digraphs.

0  of them have dichromatic number at least 3. In particular, 0 of them are 3-dicritical.

-----------------------------------------
Computing 1-extensions of the candidates on 9 vertices that are potentially subtournaments of 3-dicritical tournaments (up to isomorphism).
Process: |##################################################| 100.0% Complete
Number of  1 -extensions up to isomorphism:  34
0  of them have dichromatic number at least 3. In particular, 0 of them are 3-dicritical.

-----------------------------------------
Computing 2-extensions of the candidates on 9 vertices that are potentially subtournaments of 3-dicritical tournaments (up to isomorphism).
Process: |##################################################| 100.0% Complete
Number of  2 -extensions up to isomorphism:  0
0  of them have dichromatic number at least 3. In particular, 0 of them are 3-dicritical.
\end{lstlisting}

\subsection{The proof of Lemma~\ref{finale}}\label{b4}

We here give the code used in the proof of Lemma~\ref{finale}.

\begin{lstlisting}[language=python]
load("tools.sage")

def possible_completions(graph_to_complete, nb_vertices, list_forbidden_subdigraphs, list_forbidden_induced_subdigraphs, progress=0):
    if(progress == nb_vertices-1): 
        return [graph_to_complete]
    else:
        result = []
        for i in range(3):
            #we make a copy of the graph_to_complete
            new_D = DiGraph(nb_vertices)
            new_D.add_edges(graph_to_complete.edges())
            
            #we consider every possible orientation between the vertices (nb_vertices-1) and (progress)
            if(i==0):
                new_D.add_edge(nb_vertices-1, progress)
            elif(i==1):
                new_D.add_edge(progress, nb_vertices-1)
            else:
                new_D.add_edge(nb_vertices-1, progress)
                new_D.add_edge(progress, nb_vertices-1)
                
            #for each of the 3 possible orientations, we check whether the obtained digraph is already an obstruction. If it is not, we compute all possible completions recursively
            if(can_be_subgraph_of_3_dicritical(new_D,list_forbidden_subdigraphs, list_forbidden_induced_subdigraphs)):
                result.extend(possible_completions(new_D, nb_vertices, list_forbidden_subdigraphs, list_forbidden_induced_subdigraphs, progress+1))
        return result

for tt in range(1,8):
    transitive_tournament = DiGraph(tt)
    for i in range(tt):
        for j in range(i):
            transitive_tournament.add_edge(j,i)

    next_transitive_tournament = DiGraph(tt+1)
    for i in range(tt+1):
        for j in range(i):
            next_transitive_tournament.add_edge(j,i)

    list_forbidden_subdigraphs = [S4, K2_K2, O4, O5, K2_C3, C3_K2, C3_C3, F]
    list_forbidden_induced_subdigraphs = []
    if(tt<7):
        list_forbidden_induced_subdigraphs = [next_transitive_tournament]
    else:
        list_forbidden_induced_subdigraphs = reversed_TT8

    print("\n--------------------------------------------------------\n")
    print("Generating all 3-dicritical semi-complete digraphs with maximum acyclic induced subdigraph of size exactly " + str(tt) + ".")
    
    n=tt+1
    candidates = [transitive_tournament]
    next_candidates = []
    while(len(candidates)>0):
        print("\nComputing candidates on "+str(n)+" vertices.")
        printProgressBar(0, len(candidates))
        for i in range(len(candidates)):
            old_D = candidates[i]
            new_D = DiGraph(n)
            new_D.add_edges(old_D.edges())
            all_possible_completions_new_D = possible_completions(new_D, n, list_forbidden_subdigraphs, list_forbidden_induced_subdigraphs)
            for candidate in all_possible_completions_new_D:
                check = True
                for D in next_candidates:
                    check = not D.is_isomorphic(candidate)
                    if(not check):
                        break
                if(check):
                    next_candidates.append(candidate)
            printProgressBar(i + 1, len(candidates))

        #check the candidates that are actually 3-dicritical.
        print("We found",len(next_candidates),"candidates on "+str(n)+" vertices.")
        dicriticals = []
        for D in next_candidates:
            if(not is_two_dicolourable(D)):
                dicriticals.append(D)
        print(len(dicriticals), " of them are actually 3-dicritical.\n")
        for D in dicriticals:
            print("adjacency matrix of a 3-dicritical digraph that we found:")
            print(D.adjacency_matrix())

        candidates = next_candidates
        next_candidates = []
        n+=1
\end{lstlisting}

Running this code produces the following output after roughly 2 hours of execution on a standard desktop computer:

\begin{lstlisting}
--------------------------------------------------------

Generating all 3-dicritical semi-complete digraphs with maximum acyclic induced subdigraph of size exactly 1.

Computing candidates on 2 vertices.
Process: |##################################################| 100.0% Complete
We found 1 candidates on 2 vertices.
0  of them are actually 3-dicritical.


Computing candidates on 3 vertices.
Process: |##################################################| 100.0% Complete
We found 1 candidates on 3 vertices.
1  of them are actually 3-dicritical.

adjacency matrix of a 3-dicritical digraph that we found:
[0 1 1]
[1 0 1]
[1 1 0]

Computing candidates on 4 vertices.
Process: |##################################################| 100.0% Complete
We found 0 candidates on 4 vertices.
0  of them are actually 3-dicritical.


--------------------------------------------------------

Generating all 3-dicritical semi-complete digraphs with maximum acyclic induced subdigraph of size exactly 2.

Computing candidates on 3 vertices.
Process: |##################################################| 100.0% Complete
We found 5 candidates on 3 vertices.
0  of them are actually 3-dicritical.


Computing candidates on 4 vertices.
Process: |##################################################| 100.0% Complete
We found 5 candidates on 4 vertices.
0  of them are actually 3-dicritical.


Computing candidates on 5 vertices.
Process: |##################################################| 100.0% Complete
We found 0 candidates on 5 vertices.
0  of them are actually 3-dicritical.


--------------------------------------------------------

Generating all 3-dicritical semi-complete digraphs with maximum acyclic induced subdigraph of size exactly 3.

Computing candidates on 4 vertices.
Process: |##################################################| 100.0% Complete
We found 13 candidates on 4 vertices.
0  of them are actually 3-dicritical.


Computing candidates on 5 vertices.
Process: |##################################################| 100.0% Complete
We found 37 candidates on 5 vertices.
1  of them are actually 3-dicritical.

adjacency matrix of a 3-dicritical digraph that we found:
[0 1 1 0 0]
[0 0 1 0 1]
[0 0 0 1 1]
[1 1 0 0 1]
[1 0 1 1 0]

Computing candidates on 6 vertices.
Process: |##################################################| 100.0% Complete
We found 8 candidates on 6 vertices.
0  of them are actually 3-dicritical.


Computing candidates on 7 vertices.
Process: |##################################################| 100.0% Complete
We found 1 candidates on 7 vertices.
1  of them are actually 3-dicritical.

adjacency matrix of a 3-dicritical digraph that we found:
[0 1 1 0 0 0 1]
[0 0 1 0 1 1 0]
[0 0 0 1 0 1 1]
[1 1 0 0 0 1 0]
[1 0 1 1 0 0 0]
[1 0 0 0 1 0 1]
[0 1 0 1 1 0 0]

Computing candidates on 8 vertices.
Process: |##################################################| 100.0% Complete
We found 0 candidates on 8 vertices.
0  of them are actually 3-dicritical.


--------------------------------------------------------

Generating all 3-dicritical semi-complete digraphs with maximum acyclic induced subdigraph of size exactly 4.

Computing candidates on 5 vertices.
Process: |##################################################| 100.0% Complete
We found 27 candidates on 5 vertices.
0  of them are actually 3-dicritical.


Computing candidates on 6 vertices.
Process: |##################################################| 100.0% Complete
We found 116 candidates on 6 vertices.
0  of them are actually 3-dicritical.


Computing candidates on 7 vertices.
Process: |##################################################| 100.0% Complete
We found 10 candidates on 7 vertices.
0  of them are actually 3-dicritical.


Computing candidates on 8 vertices.
Process: |##################################################| 100.0% Complete
We found 0 candidates on 8 vertices.
0  of them are actually 3-dicritical.


--------------------------------------------------------

Generating all 3-dicritical semi-complete digraphs with maximum acyclic induced subdigraph of size exactly 5.

Computing candidates on 6 vertices.
Process: |##################################################| 100.0% Complete
We found 49 candidates on 6 vertices.
0  of them are actually 3-dicritical.


Computing candidates on 7 vertices.
Process: |##################################################| 100.0% Complete
We found 266 candidates on 7 vertices.
0  of them are actually 3-dicritical.


Computing candidates on 8 vertices.
Process: |##################################################| 100.0% Complete
We found 20 candidates on 8 vertices.
0  of them are actually 3-dicritical.


Computing candidates on 9 vertices.
Process: |##################################################| 100.0% Complete
We found 0 candidates on 9 vertices.
0  of them are actually 3-dicritical.


--------------------------------------------------------

Generating all 3-dicritical semi-complete digraphs with maximum acyclic induced subdigraph of size exactly 6.

Computing candidates on 7 vertices.
Process: |##################################################| 100.0% Complete
We found 80 candidates on 7 vertices.
0  of them are actually 3-dicritical.


Computing candidates on 8 vertices.
Process: |##################################################| 100.0% Complete
We found 500 candidates on 8 vertices.
0  of them are actually 3-dicritical.


Computing candidates on 9 vertices.
Process: |##################################################| 100.0% Complete
We found 39 candidates on 9 vertices.
0  of them are actually 3-dicritical.


Computing candidates on 10 vertices.
Process: |##################################################| 100.0% Complete
We found 0 candidates on 10 vertices.
0  of them are actually 3-dicritical.


--------------------------------------------------------

Generating all 3-dicritical semi-complete digraphs with maximum acyclic induced subdigraph of size exactly 7.

Computing candidates on 8 vertices.
Process: |##################################################| 100.0% Complete
We found 110 candidates on 8 vertices.
0  of them are actually 3-dicritical.


Computing candidates on 9 vertices.
Process: |##################################################| 100.0% Complete
We found 459 candidates on 9 vertices.
0  of them are actually 3-dicritical.


Computing candidates on 10 vertices.
Process: |##################################################| 100.0% Complete
We found 16 candidates on 10 vertices.
0  of them are actually 3-dicritical.


Computing candidates on 11 vertices.
Process: |##################################################| 100.0% Complete
We found 0 candidates on 11 vertices.
0  of them are actually 3-dicritical.
\end{lstlisting}

The adjacency matrices in the output are exactly those of the digraphs $\bid{K_3}$, $\mathcal{H}_5$ and $\mathcal{P}_7$.

\end{document}